\documentclass{notices}

\usepackage{array,url}
\usepackage{verbatim}
\usepackage{a4wide,tikz,amsmath,amssymb,amscd,xypic,amsthm,txfonts,dsfont}

\usepackage{mathtools}

\numberwithin{equation}{section}
\numberwithin{table}{section}
\numberwithin{figure}{section}
\makeatletter
\renewcommand\hline{%
  \noalign{\ifnum0=`}\fi\@ifnextchar[{\@varihline}
                                     {\@varihline[\arrayrulewidth]}}
\def\@varihline[#1]{\hrule \@height #1 \futurelet
   \reserved@a\@xhline}
\makeatother

\newcommand{\nc}{\newcommand}
\nc{\mline}{\hline[1pt]}
\nc{\acts}{\!\! \circlearrowleft}
\nc{\ccol}[1]{\multicolumn{1}{c|}{#1}}
\nc{\lccol}[1]{\multicolumn{1}{|c|}{#1}}
\nc{\gr}[1]{{\color{dg} #1}}
\nc{\re}[1]{{\color{red} #1}}
\nc{\ora}[1]{{\color{orange} #1}}
\nc{\bl}[1]{{{\color{blue}#1}}}
\nc{\ch}{\bl{\checkmark}}

\nc{\bfu}{{\bf u}}
\nc{\bfm}{{\bf m}}

\nc{\Real}{\mbox{Re}}
\nc{\End}{\mathop{\rm End}}
\nc{\Q}{{\mathbb{Q}}}
\nc{\R}{{\mathbb{R}}}
\nc{\bbP}{{\mathbb{P}}}
\nc{\C}{{\mathbb{C}}}
\nc{\G}{{\mathbb{G}}}
\nc{\A}{{\mathbb{A}}}
\nc{\I}{{\mathbb{I}}}
\nc{\M}{{\mathbb{M}}}
\nc{\Z}{{\mathbb{Z}}}
\nc{\F}{{\mathbb{F}}}
\nc{\calH}{\mathcal{H}}
\nc{\calO}{\mathcal{O}}
\nc{\calB}{\mathcal{B}}
\nc{\cM}{\mathcal{M}}
\nc{\cV}{\mathcal{V}}
\nc{\fr}{\mbox{fr}}
\nc{\ord}{\mbox{ord}}
\nc{\cX}{\mathcal{X}}
\nc{\gal}{\mathrm{Gal}}
\nc{\Cl}{\mathrm{Cl}}
\nc{\Fr}{\mathrm{Fr}}
\nc{\cl}{\mathrm{cl}}
\nc{\Gal}{\mathrm{Gal}}
\nc{\OK}{\mathcal{O}_K}
\nc{\cK}{\mathcal{K}}
\nc{\OKp}{\mathcal{O}_{K'}}
\nc{\Kgal}{K^g}
\nc{\grd}{\textrm{grd}}
\nc{\lcm}{\textrm{lcm}}
\nc{\cmt}[1]{}
\nc{\kpg}{\mathcal{K}_{G,p}}
\nc{\pco}{p_G; \tame{p_G}}
\nc{\sd}{{:}}
\nc{\pstar}{\hat{p}}
\nc{\cP}{\mathcal{P}}
\nc{\cd}{\! \cdot \!}
\nc{\z}{\ }
\nc{\qqud}{\strut\qquad}
\nc{\T}{\mathbb{T}}
\nc{\cmmt}[1]{}
\nc{\para}[1]{\vspace{-.17in} \paragraph{#1}}
\nc{\ccdot}{\! \cdot \!}
\nc{\w}{\hspace{5pt}}
\nc{\Sp}{\mbox{Sp}}
\nc{\Or}{\mbox{O}}
\nc{\CSp}{\mbox{CSp}}
\nc{\COr}{\mbox{CO}}
\nc{\fa}{Q}
\nc{\faa}{q}

\DeclareMathOperator{\Tr}{Tr}

\DeclareMathOperator{\vol}{vol}

\DeclareMathOperator{\tame}{tame}

\def\vdts{\raisebox{2pt}[0pt]{\vbox{\baselineskip2pt \lineskiplimit0pt
    \kern6pt\hbox{.}\hbox{.}\hbox{.}}}}
\nc{\vdd}{-4pt}

\allowdisplaybreaks

\date{}

\title{Hypergeometric Motives}

\author{ David P.\ Roberts \\ Fernando Rodriguez Villegas}

\begin{document}

\maketitle

\section{Introduction}
\label{intro}

It must have been frustrating in the early days of calculus that an integral like 
\begin{equation}
\label{legendre-integral}
F(t)=\frac{1}{\pi} \int_0^1\frac1{\sqrt{x(1-x)(1-tx)}}\,dx
\end{equation}
 appeared not to be
expressible in terms of known functions. 
This type of integral arises in computing the movement of the ideal pendulum 
or the length of an arc of an ellipse for example; they have remained
relevant and are connected to a great deal of the mathematics of the
last 200 years.

Indeed $F$ is not an elementary function. Its Maclaurin expansion 
\begin{equation}
\label{legendre-series}
F(t)=\sum_{k =0}^\infty \binom {2k}k^2\, \left(\frac t{16}\right)^k
\end{equation}
is an example of a hypergeometric series.  It satisfies a
linear differential equation of order two of the type brilliantly
analyzed by Riemann. As mentioned by Katz~\cite[p.3]{Katz-RLS}, Riemann
was lucky.  His analysis only works because any rank two differential equation on
$\bbP^1(\C) - \{0,1,\infty\}$
is {\it rigid} in the sense that the local behavior of solutions around
the missing points uniquely determines their global behavior. 

Taking a more geometric perspective, \eqref{legendre-integral}
is presenting the function $\pi F$ as a period of the 
family of elliptic curves defined by
\begin{equation}
\label{legendre-equation}
E_t: \quad y^2= x(1-x)(x-t).
\end{equation}
This fact implies as well
 that $F$
satisfies an order two linear differential equation, ultimately
because $H^1(E_t,\Q)$ is two-dimensional.  \let\thefootnote\relax\footnote{David P. Roberts is a professor at the University of Minnesota Morris.  His e-mail address is {\tt roberts@morris.umn.edu} and his research is supported by  grant DMS-1601350 from the National Science Foundation.}
\footnote{Fernando Rodriguez Villegas is a senior research scientist at the Abdus Salam International Centre for Theoretical Physics.  His email address is {\tt villegas@ictp.it}.}

Shifting now to more arithmetic topics, 
if we fix a rational number $t\neq 0,1$
then for almost all primes $p$ the number $a_p$ defined by
\begin{equation}
\label{legendre-count}
|E_t(\F_p)|=p+1-a_p
\end{equation}
is of fundamental importance.  With these $a_p$ as the main ingredients, one builds
an $L$-function 
\begin{equation}
\label{legendre-L} L(E_t,s) = \sum_{n=1}^\infty \frac{a_n}{n^s}.
\end{equation}
Much of the importance of the $a_p$ is seen through this 
$L$-function.  For example, the famous Birch-Swinnerton-Dyer
conjecture says that the group $E_t(\Q)$ modulo its torsion
is isomorphic to $\Z^r$, where $r$ is the order of vanishing 
of $L(E_t,s)$ at $s=1$.  A critical advance is the result of Wiles {\em et al.}\ that
the function 
\begin{equation}
\label{legendre-mod}
f(z) = \sum_{n=1}^\infty a_n e^{2 \pi i z} 
\end{equation}
on the upper half plane is a modular form.  In particular, this result
implies that $L(E_t,s)$ is at least well-defined at $s=1$. 

The equations displayed so far represent a standard general paradigm in 
arithmetic geometry.  One can start with any variety $X$ over $\Q$, not just the 
varieties \eqref{legendre-equation}.  There are fully developed
theories of periods and point counts, and in principle one can
produce analogs of the period formulas \eqref{legendre-integral}-\eqref{legendre-series} 
and the point count formula \eqref{legendre-count}. Interacting now with deep but widely-believed
conjectures, one can break the cohomology of $X$ into 
irreducible motives, study $L$-functions
like \eqref{legendre-L}, and try to find corresponding
automorphic forms like \eqref{legendre-mod}. 

This survey is an informal invitation to {\em hypergeometric
  motives}, hereafter abbreviated HGMs; see Section~\ref{HGMs} for their
  definition. We write them as $H(\fa,t)$, with a rational function
$\fa \in \Q(T)$ satisfying certain conditions being the family
parameter and $t \in \Q-\{0,1\}$ the specialization parameter.  The
introductory family of examples is
\begin{equation} 
H \left({(T+1)^2}/{(T-1)^2},t \right) = H^1(E_t,\Q).
\end{equation}
  Rigidity makes HGMs much more tractable than general motives:
  periods, point counts, and other invariants are given by explicit
  formulas in the parameters $(\fa,t)$.  Our broader goal in this survey is to use HGMs to gain
 insight into the general theory of motives; we
 illustrate all topics with explicit examples throughout. 

Sections \ref{monodromy}-\ref{special-motive} are geometric in nature.  The main focus is 
on varieties generalizing \eqref{legendre-equation} 
and the discrete aspect of periods like \eqref{legendre-integral}-\eqref{legendre-series},
as captured in vectors of 
Hodge numbers, $h = (h^{w,0},\dots,h^{0,w})$.  
A theme here is that HGMs form quite a broad class of
irreducible motives, as very general $h$ arise.   
 Sections~\ref{sect:arithmetic}-\ref{numerical} are 
arithmetic in nature, with the focus being on
generalizations of~\eqref{legendre-count},~\eqref{legendre-mod},
and especially~\eqref{legendre-L}.   Watkins 
has written a very useful hypergeometric motives package \cite{Wat}
in {\em Magma} and throughout this article we indicate
how to use it by including small snippets of {\em Magma} code.
Together these snippets are enough to let {\em Magma}
beginners numerically compute with $L$-functions $L(H(\fa,t),s)$ using the
free online {\em Magma} calculator.

\section{Hypergeometric functions}
\label{monodromy}
We begin by generalizing \eqref{legendre-integral}-\eqref{legendre-series} and explaining how
this generalization leads to family parameters.  

\para{Integrals and series.} Let $\alpha = (\alpha_1, \dots, \alpha_n)$, 
and $\beta = (\beta_1, \dots, \beta_n)$ be vectors of complex numbers
with $\Real(\beta_j)>\Real(\alpha_j)>0$ and $\beta_n=1$.  For $|t| < 1$ define, making use of the
standard Gamma function
$\Gamma(s) = \int_0^\infty e^{-x} x^{s-1} dx$, 
\begin{eqnarray}
\label{euler-integral}
\lefteqn{F(\alpha,\beta,t)=  \prod_{i=1}^n\frac{\Gamma(\beta_i)}{\Gamma(\alpha_i)\Gamma(\beta_i-\alpha_i)} \cdot} \\
&& \int_0^1\cdots\int_0^1  \frac{\prod_{i=1}^{n-1} \left( x_i^{\alpha_i-1}(1-x_i)^{\beta_i-\alpha_i-1}  dx_i \right)}{ (1-tx_1\cdots x_{n-1})^{\alpha_n}}. \nonumber
\end{eqnarray}
Via $\Gamma(1)=1$ and $\Gamma(1/2) = \sqrt{\pi}$,  \eqref{legendre-integral} is the special case $\alpha = (1/2,1/2)$ and $\beta=(1,1)$.

Expand the denominator of the integrand of \eqref{euler-integral} via the binomial theorem and use
Euler's beta integral to evaluate the individual terms.  Written in terms 
of Pochhammer symbols $(u)_k 
= u(+1) \cdots (u+k-1)$, the result is
\begin{equation}
\label{hyperg-series-defn}
F (\alpha,\beta,t) = {}_{n}F_{n-1} (\alpha,\beta,t):=
   \sum_{k = 0}^\infty \frac{(\alpha_1)_k \cdots (\alpha_n)_k}{(\beta_1)_k\dots
  (\beta_n)_k} t^k.
\end{equation}
In other words, the integral \eqref{euler-integral} is an alternative
definition of the standard hypergeometric power series \eqref{hyperg-series-defn}.
The case $\alpha = (1/2,1/2)$ and $\beta=(1,1)$ simplifies to \eqref{legendre-series}.

\para{Monodromy.} An excellent general reference for hypergeometric functions is \cite{BH},
and we now give a summary sufficient for this survey.   The function
$F(\alpha,\beta,t)$ is in the kernel of an $n^{\rm th}$ order 
differential operator $D(\alpha,\beta)$ with singularities 
only at $0$, $1$ and $\infty$.   This means in particular 
that $F(\alpha,\beta,t)$, initially defined on the unit disk, 
extends to a ``multivalued function'' on the thrice-punctured
projective line $\bbP^1(\C)-\{0,1,\infty) = \C-\{0,1\}$.   With respect to a given
basis, this multivaluedness is codified by a representation $\rho$ of the fundamental group
$\pi_1(\C-\{0,1\},1/2)$ into $GL_n(\C)$.   The 
fundamental group is free on $g_0$ and $g_1$, 
with these elements coming from counterclockwise 
circular paths of radius $1/2$ about $0$ and $1$ 
respectively.  To emphasize the equal status of $\infty$ and $0$, it is better
to present this group as generated by $g_\infty$, $g_1$, 
and $g_0$, subject to the relation 
$g_\infty g_1 g_0=1$.   The assumption 
$\beta_n=1$ was only imposed to present
the classical viewpoint 
\eqref{euler-integral}-\eqref{hyperg-series-defn} cleanly;
we henceforth drop it.

A useful fact due to Levelt is the explicit description of the matrices 
$h_\tau=\rho(g_\tau) \in GL_n(\C)$ with respect to a certain well-chosen basis.  
Define polynomials
\begin{align*}
\faa_\infty & := (T-e^{2\pi i \alpha_1})\cdots(T-e^{2\pi i\alpha_n}), \\
\faa_0 &:=(T-e^{-2\pi i \beta_1})\cdots(T-e^{-2\pi i  \beta_n}).
\end{align*}
Then $h_\infty$ and $h_0$ 
are companion matrices of $\faa_\infty$ and $\faa_0$, while $h_1$
is determined by $h_\infty h_1 h_0=I$.   The
matrix $h_1$ differs minimally from 
the identity in that $h_1-I$ has rank
$1$.   We will henceforth consider only cases where
no $\alpha_i-\beta_j$ is an integer.  This ensures
that the $h_\tau$ generate an irreducible 
subgroup $\Gamma$ of $GL_n(\C)$.  Moreover the 
representation is rigid, in the following
sense: suppose $h_\infty'$, $h'_1$, and $h_0'$
are conjugate to $h_\infty$, $h_1$, and $h_0$
respectively.  Then there is a single 
matrix $c$ such that $h_\tau = ch_\tau'c^{-1}$ 
for all three $\tau$.  

\para{Family parameters.}
    The parameters $(\alpha,\beta)$ 
contain information which is irrelevant for the
sequel.   First, the individual $\alpha_j$ and
$\beta_j$ are important only modulo integers.
Second, the orderings of the $\alpha_j$ and $\beta_j$
  do not matter.  To remove these irrelevancies, we will regard the
degree $n$ rational function $\fa = \faa_\infty/\faa_0$ as the primary
index in the sequel, calling it the {\em family parameter}.  A bonus
of this shift in emphasis is that an important field $E \subset \C$ is
made evident, the field generated by the coefficients of $\faa_\infty$
and $\faa_0$.  By construction, all three $h_\tau$ lie in $GL_n(E)$.

    The cases which naturally have 
underlying motives are exactly the ones
with all $\alpha_j$ and $\beta_j$ rational,
so that $E$ is some cyclotomic field.  In
this survey we will substantially simplify by
restricting to cases with $E = \Q$.
With this simplification, there are two 
natural ways to present $\fa$ as follows.
Write $\Psi_m=T^m-1$ and consider 
its factorization into irreducible 
polynomials, $\Psi_m = \prod_{d|m} \Phi_d$. 
So the factors are cyclotomic polynomials
$\Phi_d = \prod_{j \in (\Z/d)^\times} (T - e^{2 \pi i j/d})$
and thus have degree the 
Euler totient $\phi(d) = |(\Z/d)^\times|$.   

In our introductory example, the ways are
\begin{equation}
\label{qphipsi}
\fa   \;\;\; = \;\;\; \frac{(T^2-1)^2}{(T-1)^4} = \frac{\Psi^2_2}{\Psi_1^4} \;\;\; =  \;\;\; \frac{(T+1)^2}{(T-1)^2} = \frac{\Phi^2_2}{\Phi_1^2}.
\end{equation}
In general, the second way is just the canonical factorization into irreducibles, while the first way is the unique ``unreduction'' to
products of $\Psi_m$ in which no factor appears in both a numerator and denominator.

To enter a family parameter $\fa$ into {\em Magma}, one
can use either of these two ways, as in 
the equivalent commands
\begin{eqnarray}
&&\mbox{\tt Q:=HypergeometricData([*-2,-2,1,1,1,1*]);} \nonumber  \\ 
\label{defineq} &&\mbox{\tt Q:=HypergeometricData([1,1],[2,2]);} 
\end{eqnarray}
In the first method,  one
inputs just the {\em gamma vector} $\gamma = [\gamma_1,\dots,\gamma_l]$ formed by subscripts on the $\Psi$'s, 
using signs to distinguish between numerator and denominator.
In the second method, one inputs just the subscripts of the 
denominator and then numerator 
$\Phi$'s, 
these being called the {\em cyclotomic parameters}.  
When working with underlying varieties, the gamma vectors
are so useful that we often simply write $H(\gamma,t)$ rather than $H(\fa,t)$. After the transition to motives,
the cyclotomic presentation is generally more convenient.  To simplify
slightly, we henceforth require that $\gcd(\gamma_1,\dots,\gamma_l) = 1$.

Note that initialization commands like \eqref{defineq} do nothing by themselves;
in this survey {\em Magma} will first start returning useful information in 
Sections~\ref{sect:hodge-numbers} and \ref{frobenius}.   Note also that
{\em Magma} requires a semicolon at the end of all commands,
as in \eqref{defineq}.  We often omit these semicolons in the sequel.

\para{Orthogonal vs.\ symplectic.}
The number $\fa(0)=\det(h_1)$ is either $-1$ or $1$ under our
restriction $E=\Q$.  This dichotomy is strongly felt throughout this
survey.  It can also be expressed in terms of the fundamental
bilinear form $\left( \cdot, \cdot \right)$ on $\Q^n$ preserved by the
monodromy group $\Gamma = \langle h_\infty,h_0 \rangle$; 
see~\cite[\S4]{BH},~\cite[\S3.5]{RV-bez}.
In the {\em orthogonal case}, $h_1$ is conjugate to
$\mbox{diag}(-1,1,\dots,1)$ and $\left( \cdot,\cdot \right)$ is
symmetric.  In the {\em symplectic case}, $h_1$ is conjugate to
$\binom{1 \; 1}{0 \; 1} \oplus \mbox{diag}(1,\dots,1)$, and
$\left( \cdot,\cdot \right)$ is antisymmetric.

\section{Source varieties} 
\label{varieties} We now describe varieties which give rise 
to hypergeometric motives. 

\para{Euler varieties.} We have already generalized \eqref{legendre-integral}
to \eqref{euler-integral} and \eqref{legendre-series} to 
\eqref{hyperg-series-defn}.    Assuming briefly $\beta_n=1$ again, 
a natural generalization of \eqref{legendre-equation} is
to 
\begin{equation}
\label{eulercover}
y^m=\prod_{j=1}^{n-1}x_j^{a_j}(1-x_j)^{b_j}(1-tx_1\cdots x_{n-1})^{a_n}.
\end{equation}
Here $m$ is the least common denominator of the $\alpha_j$ and $\beta_j$,
and the exponents are integers $0\leq a_j,b_j<m$ such that 
$$
a_j\equiv -m\alpha_j\bmod m, \qquad b_j\equiv
m(\alpha_j-\beta_j)\bmod m,
$$ 
for $j=1,\ldots,n-1$ and $a_n\equiv m \alpha_n\bmod m$.  
The equations \eqref{euler-integral}-\eqref{hyperg-series-defn} show that a 
specified scalar multiple of ${}_nF_{n-1}(\alpha,\beta,t)$ arises as a period of 
this variety.  However the varieties \eqref{eulercover} depend on how
the parameters are paired: 
$(\alpha_1,\beta_1)$, \dots, $(\alpha_n,\beta_n)$.
This dependence complicates the arithmetic 
of these varieties, so we will use an alternative
collection of varieties to define 
hypergeometric motives.  

{\para{Canonical varieties.} The alternative varieties 
appear under the term ``circuits'' in \cite{GKZ} and 
are studied at greater length in \cite{BCM}. 
For a gamma vector $\gamma$ and a complex number $t$, define  
$X^{\rm bcm}_{\gamma,t} \subset \bbP^{l-1}$ by 
two homogeneous equations,
\begin{align}
\label{bcm}
\sum_{j=1}^{l} y_j & = 0, &
\prod_{\gamma_j>0}y_j^{\gamma_j}= u \prod_{\gamma_i<0}y_j^{-\gamma_j}.
\end{align}
Here and in the sequel, we systematically use the abbreviation $u = t \prod_j \gamma_j^{\gamma_j}$.  
The canonical variety is by definition the open subvariety $X_{\gamma,t}$ 
on which all the homogeneous coordinates $y_j$ are nonzero.   
The point 
$(\gamma_1: \cdots : \gamma_l)$ is
an ordinary double point on $X_{\gamma,1}$ and
otherwise all the $X_{\gamma,t}$ are smooth.   Because
of this double point, we exclude the case $t=1$ from consideration
until Section~\ref{special-motive}.

\para{Toric models.} From a dimension-count viewpoint, the BCM
equations \eqref{bcm} for canonical varieties are inefficient.  They
start with the $l = \kappa+3$ variables $y_i$ and use two equations
and projectivization to get the desired $\kappa$-dimensional variety
$X_{\gamma,t}$.  The toric models from \cite{GKZ} start
instead with $d=\kappa+1$ variables $x_i$ and present $X_{\gamma,t}$
by just one equation.

To obtain a toric model from a gamma vector
$\gamma$, one proceeds as illustrated by Table~\ref{toroidalchart}. 
First, for each new variable $x_i$ choose a row vector $m_{i*}$ in $\Z^l$ which
is orthogonal to the given $l$-vector $\gamma$.  These row vectors are
required to be such that $\Z^l/\langle m_{i*} \rangle$ is torsion-free.  
Second, choose a row vector $k \in \Z^l$ which satisfies $\gamma \cdot k = 1$.  
The toric model is then 
\begin{equation}
\label{toroidal}
\sum_{i=1}^l u^{k_i} \prod_{j=1}^d x_i^{m_{ij}}= 0.
\end{equation}
So in the example of Table~\ref{toroidalchart}, the  
resulting equation is 
\begin{equation}
\label{toroidalexample}
x_1^2 + u x_1 x_2^2 + u + x_1^3 x_2= 0,
\end{equation}
with $u=-2^63^3 t/5^5$.    In general, the variety $X_{\gamma,t}$ is the subvariety of 
the torus $\G_m^d$ given by the equation \eqref{toroidal}.
 
 \begin{table}[htb]
\[
{\renewcommand{\arraycolsep}{3pt}
\begin{array}{r|cccc|}
\cline{2-5}
&  \gamma_1 & \gamma_2 & \gamma_3 & \gamma_4 \\
 \cline{2-5}
x_1:  \;\;\; &  m_{11} & m_{12} & m_{13} & m_{14} \\
x_2:  \;\;\; & m_{21} & m_{22} & m_{23} & m_{24} \\
 \cline{2-5}
u: \;\;\; &  k_1 & k_2 & k_3 & k_4 \\
 \cline{2-5}
\end{array} 
=
\begin{array}{|cccc|}
\hline
  -5 & - 2 & 3 & 4 \\
 \hline
  2 & 1 & 0 & 3 \\
  0 & 2 & 0 & 1 \\
 \hline
 0 & 1 & 1 & 0 \\
 \hline
\end{array} }
\]
\caption{\label{toroidalchart} Derivation of the equation \eqref{toroidalexample} for $X_{[-5,-2,3,4],t}$}
\end{table}

The relation between the BCM equation for $X_{\gamma,t}$ and 
a toric model for $X_{\gamma,t}$ is very simple:
\begin{equation}
\label{parameterization}
y_j =  u^{k_j} \prod_{i=1}^d x_i^{m_{ij}}.
\end{equation} 
When one uses \eqref{parameterization} to write \eqref{bcm} in terms of the $x_i$,
the second equation is identically satisfied while the first becomes \eqref{toroidal}.  
Conversely, any point $(y_1: \dots : y_l)$ comes from a unique $(x_1,\dots,x_d)$ 
because of the torsion-free condition.

 \para{Polytopes.} A toric model gives a polytope $\Delta \subset \R^d$ 
 which is an aid to understanding  
 the $X_{\gamma,t}$.   The case $d=2$ is readily
 visualizable and Figure~\ref{toroidalpict} 
 continues our example.  In general, one 
 interprets the column vectors $m_{*j}$ of
 the chosen matrix as points in $\Z^d$ and 
 $\Delta$ is their convex hull.  Let 
 $\Delta_j$ be the convex hull of
 all the points except the $j^{\rm th}$ 
 one.  
Normalize volume so that the standard $d$-dimensional simplex has
  volume~$1$, and thus $[0,1]^d$ has volume $d!$.
 Then the volume
 of $\Delta_j$ is $|\gamma_j|$.   
 The $\Delta_j$ with $\gamma_j>0$ form
 one triangulation of $\Delta$, while
 the $\Delta_j$ with $\gamma_j<0$ form
 another.   The total volume of $\Delta$ is
 the important number $\mbox{vol}(\gamma)= \frac{1}{2} \sum_{j=1}^l |\gamma_j|$
 
 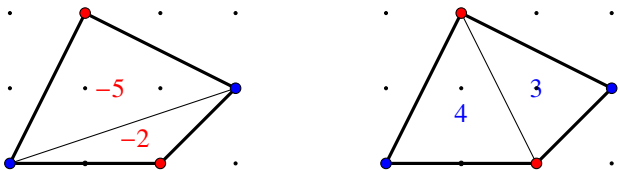
\begin{figure}[htb]
\begin{center}
\begin{tikzpicture}
\draw[very thick] (0,0) -- (1,2);
\draw[very thick] (0,0) -- (2,0);
\draw[very thick] (1,2) -- (3,1);
\draw[very thick] (3,1) -- (2,0);
\draw[thin] (0,0) -- (3,1);
\draw[red] (1.66,.33) node{$-2$};
\draw[red] (1.33,1)  node{$-5$};
\draw[fill=blue] (0,0) circle (.07);
\draw[fill=red] (1,2) circle (.07);
\draw[fill=blue] (3,1) circle (.07);
\draw[fill=red] (2,0) circle (.07);
\draw[fill=black] (0,1) circle (.02);
\draw[fill=black] (0,2) circle (.02);
\draw[fill=black] (1,0) circle (.03);
\draw[fill=black] (1,1) circle (.02);
\draw[fill=black] (2,1) circle (.02);
\draw[fill=black] (2,2) circle (.02);
\draw[fill=black] (3,0) circle (.02);
\draw[fill=black] (3,2) circle (.02);

\draw[very thick] (5,0) -- (6,2);
\draw[very thick] (5,0) -- (7,0);
\draw[very thick] (6,2) -- (8,1);
\draw[very thick] (8,1) -- (7,0);
\draw[thin] (7,0) -- (6,2);
\draw[blue] (6,.67) node{$4$};
\draw[blue] (7,1) node{$3$};

\draw[fill=blue] (5,0) circle (.07);
\draw[fill=red] (6,2) circle (.07);
\draw[fill=blue] (8,1) circle (.07);
\draw[fill=red] (7,0) circle (.07);
\draw[fill=black] (5,1) circle (.02);
\draw[fill=black] (5,2) circle (.02);
\draw[fill=black] (6,0) circle (.03);
\draw[fill=black] (6,1) circle (.02);
\draw[fill=black] (7,1) circle (.02);
\draw[fill=black] (7,2) circle (.02);
\draw[fill=black] (8,0) circle (.02);
\draw[fill=black] (8,2) circle (.02);

\end{tikzpicture}
\end{center}
\caption{\label{toroidalpict} The triangulations $\Delta = \Delta_1\cup \Delta_2$ and $\Delta=\Delta_3 \cup \Delta_4$ of 
the polytope $\Delta$ for the family with $\gamma = [-5,-2,3,4]$.  The points are at the column vectors $\binom{m_{1j}}{m_{2j}}$ of 
Table~\ref{toroidalchart}, and $\gamma_j$ is printed in the opposite triangle.}
\end{figure}

The common topology of the $X_{\gamma,t}$ with $t \neq 1$
  is reflected in the combinatorics of $\Delta$.  In the 
  case of $d=2$, the genus $g$ of $X_{\gamma,t}$ is the
  number of lattice points on the interior, while the number
  of punctures $k$ is the number of lattice points on the 
  boundary.  Pick's theorem then says that the
  Euler characteristic $\chi=2-2g-k$ of $X_{\gamma,t}$ 
  is $-\mbox{vol}(\gamma)$.  In the
  example of Figure~\ref{toroidalpict}, $(g,k,\chi) = (2,5,-7)$.  
  For larger ambient dimension $d$, the situation is of course much more
  complicated, but always $\chi = (-1)^{d-1} \mbox{vol}(\gamma)$.

  \para{Compactifications.} In algebraic geometry, one normally wants
  to compactify a given open variety such as $X_{\gamma,t}$ and there
  are typically many natural ways of doing it.  

  We already saw the compactification $X_{\gamma,t}^{\rm bcm}$.  It is
  a hypersurface of degree $\mbox{vol}(\gamma)$ in the projective
  space $\bbP^{d}$ defined by the first equation of \eqref{bcm}.  On
  the other hand, for any choice of matrix $m$ with all entries
  nonnegative, homogenization of \eqref{toroidal} gives a
  alternative compactification
  $\overline{X}_{\gamma,t} \subset \bbP^d$.
  
  In our continuing example $\gamma = [-5,-2,3,4]$, the plane curve
  $X_{\gamma,t}^{\rm bcm}$ has degree seven.  In contrast, the plane
  curve $\overline{X}_{\gamma,t}$ has degree just four, this number
  arising as the maximum column sum of the matrix $m$ in
  Table~\ref{toroidalchart}.  Smooth curves in these degrees have
  genera $15$ and $3$ respectively.  For $t \neq 1$, $X_{\gamma,t}$
  has genus $2$ so $X_{\gamma,t}^{\rm bcm}$ must have bad
  singularities while $\overline{X}_{\gamma,t}$ has just a single node.  
  
  Another compactification $X_{\gamma,t}^{\rm BCM}$ is a major focus
  of \cite{BCM}. It is typically not smooth, but only has quotient
  singularities.  These singularities are mild in
  the sense that $X_{\gamma,t}^{\rm BCM}$ looks smooth from the
  viewpoint of rational cohomology, and may be ignored when discussing
  motives as in the next section.  

\section{HGMs from cohomology}
\label{HGMs}
   Here we define HGMs and explain how their
 behavior is simpler than other similar motives.

\para{Motivic formalism.}  Let $K$ and $E$ be 
subfields of $\C$; the case of principal interest to us is $K=E=\Q$.
  Minimally modifying Grothendieck's original conditional definitions, 
Andr\'e unconditionally defined a category $\cM(K,E)$ of 
pure motives over $K$ with coefficients in $E$ \cite{And}.  
The formal structures of this category can best
be understood in terms of a huge proreductive
algebraic group $\G_K$ over $\Q$, the absolute motivic
Galois group of $K$.  Then $\cM(K,E)$ is exactly
the category of representations of $\G_K$ on
finite-dimensional $E$ vector spaces.  

When taking cohomology, we are always implicitly
working with the complex points of a variety.  
For a smooth projective
variety $X$ over $K$ and an integer $w$, 
the singular cohomology space $M=H^w(X,E)$ 
is an object of $\cM(K,E)$.  The image $G_M$
of $\G_K$ in the general linear group
of $M$ is by definition the motivic Galois group of $M$.  
The purpose of $G_M$, as the rest of this survey will
make clear, is to group-theoretically coordinate
very concrete structures on the vector space
$H^w(X,E)$.

Two copies of the multiplicative group $\G_m = GL_1$
play important roles in the formalism 
of motives.    One is a normal subgroup and the 
other a quotient: $\G_m \subset \G_K \twoheadrightarrow \G_m$.
A rank $n$ motive $M$ is said to be of weight
$w$ if the representation restricted to 
the subgroup $\G_m$ consists of $n$ copies the 
representation $r \mapsto r^w$.   The motives
$H^w(X,E)$ all have weight $w$.  
The representation of $\G_K$ on the rank one motive $E(-1) := H^2(\bbP^1,E)$ 
corresponds to the representation $t \mapsto t$
of the quotient group $\G_m$.   The motive corresponding
to the representation $t \mapsto t^j$ is denoted $E(-j)$.  
The motives $M(j) := M \otimes E(j)$ are called the 
Tate twists of $M$. 

Each category of pure motives $\cM(K,E)$ is contained 
in a larger category $\cM \cM(K,E)$ of mixed motives, 
where an irreducible motive has a canonical weight filtration 
with subquotients in $\cM(K,E)$.  We will 
mention mixed motives at several junctures, 
but our focus is sharply on pure motives.

\para{Definition of HGMs.}  Let $\gamma$ be a gamma vector
  of length $\kappa+3$ with $r$ negative entries and let
  $t \in \Q^\times-\{1\}$. The hypergeometric motive $H(\gamma,t)$ is
  defined from the cohomology of the affine variety
  $X_{\gamma,t}$~\cite{RV-Mixed}.  We start with the compactly supported middle cohomology space
  $H_c^\kappa(X_{\gamma,t},\Q)$ and begin by cutting out a subquotient $H'(\gamma,t)$ in two steps.

First, we eliminate the contribution of the ambient $d$-dimensional torus to
obtain the primitive subspace $PH_c^\kappa(X_{\gamma,t},\Q)$.
Second, we take any smooth compactification $\bar X$ of $X_{\gamma,t}$, or one with at worst 
mild singularities as mentioned above, and
consider the image $H'(\gamma,t)$ of $PH_c^\kappa(X_{\gamma,t},\Q)$ under the
natural map to $H^\kappa(\bar{X},\Q)$.  As a quotient of $PH_c^\kappa(X_{\gamma,t},\Q)$, the space
$H'(\gamma,t)$ is independent of the choice of compactification.

For example, for the compactification $X^{BCM}_{\gamma,t}$ there is a
decomposition of its middle cohomology,
$H^{\kappa}(X^{BCM}_{\gamma,t},\Q) = H'(\gamma,t) \oplus T$. It is described
at the level of point counts in \cite[Thm~1.5]{BCM}.  Here $T$ is zero
if $\kappa$ is odd and the sum of $\binom{\kappa+1}{r-1}$ copies of
$\Q(-\kappa/2)$ if $\kappa$ is even.

Finally, we define the hypergeometric motive $H(\gamma,t) \in \cM(\Q,\Q)$ as the Hodge-normalized Tate
twist $H'(\gamma,t)(j)$, as discussed in the next section.  So
$H(\gamma,t)$ has weight $w = \kappa-2j$ with $j \in \Z_{\geq 0}$
specified there.  

More conceptually, $PH_c^\kappa(X_{\gamma,t},\Q)$ is a mixed motive of
rank $\vol(\gamma)-1$ and the pure motive $H'(\gamma,t)$ is its weight
$\kappa$ quotient.  The passage from $PH_c^\kappa(X_{\gamma,t},\Q)$ to
$H'(\gamma,t)$ is closely related to the reduction of fractions as in
\eqref{qphipsi}.  In particular, $H'(\gamma,t)$ has rank
$n = \deg(\fa)$.

It is worth stressing that the full mixed motive
$PH_c^\kappa(X_{\gamma,t},\Q)$ is itself of great interest, with its
lower weight parts playing an important role in deeper studies of
hypergeometric motives.

\para{Motivic Galois groups of HGMs.} For $t \in \C^\times-\{1\}$, one likewise
gets a motive $M = H(\gamma,t) \in \cM(\Q(t),\Q)$. lf $t$ is transcendental 
then the motivic Galois group of $M$
can be cleanly expressed in terms of the monodromy
group $\Gamma$ of Section~\ref{monodromy}, as follows.   
 If $w=0$, then $\Gamma$ is finite and $G_M=\Gamma$.  
 If $w>0$, then $\Gamma$ is infinite and $G_M$ is 
 the smallest algebraic group containing both $\Gamma$
 and scalars; more explicitly, $G_M$ is the conformal
 symplectic group $\CSp_n$ in the symplectic
 case of odd $w$, and a conformal orthogonal group
 $\COr_n$ in the orthogonal case of even $w$. 
  In the case that $t$ is algebraic, including our main case that
 $t$ is rational, the same identification of $G_M$
 holds almost always.

\para{Related motives.}  The
toric model viewpoint is part of 
the program in  \cite{GKZ} to approach
algebraic geometry by emphasizing
the number $\kappa+\epsilon$ of terms 
in polynomials defining $\kappa$-dimensional 
varieties.   By scaling to normalize coefficients,
such varieties come in $(\epsilon-2)$-dimensional
families.  

  For $\epsilon=2$, the Newton polytope $\Delta$ is a simplex. 
An abelian group $A$ of order $\mbox{vol}(\Delta)$ 
and some exponent $m$ acts on the single associated 
complex variety $X$.
The essential cases here are the Fermat varieties in 
$\bbP^{\kappa+1}$, defined by 
\begin{equation}
\label{fermat}
x_1^{m} + \cdots + x_{\kappa+2}^m=0.
\end{equation}
The group $A$ comes from scaling the 
variables by $m^{\rm th}$ roots of unity
and has order $m^{\kappa+1}$ and 
exponent $m$.   Writing $K=\Q(e^{2 \pi i/m})$,
the action decomposes $H^{\kappa}(X,K)$
into one-dimensional motives in 
$\cM(K,K)$.   This setting of $\epsilon=2$ 
was the focus of several influential papers 
of Weil from around 1950, and
the rank one motives appearing 
are {\em Jacobi motives.} 

   The case $\epsilon=3$ corresponds to {\em general
hypergeometric motives} where the $\alpha_j$ 
and $\beta_j$ can be arbitrary rational numbers.   The group $A$ now
has order a divisor of $\mbox{vol}(\Delta)$ 
and some exponent $m$.   For example, 
for $m=\kappa+2 \geq 3$ one can add the term
 $u x_1 \cdots x_{m}$ to \eqref{fermat};
 then $A$ is reduced to having order 
$m^{\kappa}$ but still has exponent $m$.  
The action of $A$ again decomposes 
$H^{\kappa}(X,K)$ in $\cM(K,K)$ and the 
summands include general hypergeometric
motives.  Our torsion-free requirement for exponent matrices is 
equivalent to requiring that the column vectors 
affine span $\Z^d$; in turn, this means that our
HGMs constitute exactly the case $|A|=1$.

Much of what we are describing in this  
article both has simpler analogs for Jacobi motives and 
extends to general hypergeometric motives.  
  Indeed \cite{BH} and 
\cite{Katz-ESDE} are in the latter setting.  
However the associated 
$L$-functions correspond to 
motives that have been descended
to $\cM(\Q,\Q)$ and have rank 
$\phi(m) n$.  Because of the factor
$\phi(m)$, inclusion of these other settings
would only 
modestly increase the 
collection of computationally accessible $L$-functions.
Also the resulting motives in $\cM(\Q,\Q)$
have motivic Galois groups which
are more complicated than 
the $\CSp_n$ and $\COr_n$ arising ubiquitously
in our setting of $H(\gamma,t)$.

\section{Hodge numbers}
\label{sect:hodge-numbers}
   One of the very first things one wants to know about a
motive is its Hodge numbers.  Fortunately, this desire is easily
satisfied for HGMs by an appealing procedure.     
\para{Background.}   For a smooth projective variety $X$ over 
$K \subseteq \C$, 
there is a decomposition of complex vector spaces
$
H^w(X,\C) = \bigoplus_{p=0}^w H^{p,w-p}
$, with complex conjugation
on coefficients switching $H^{p,q}$ and $H^{q,p}$.  
 The Hodge numbers $h^{p,q} := \dim(H^{p,q})$ therefore satisfy 
Hodge symmetry $h^{p,q}=h^{q,p}$ and sum to the 
Betti number $b_w := \dim(H^w(X,\C))$.   Classical examples are given
in \eqref{bettipluri}-\eqref{projhodge} below.  

Likewise, the rank of a weight $w$ motive 
$M \in \cM(K,E)$ is decomposed into Hodge 
numbers $h^{p,w-p}$.  The decomposition has a simple 
group-theoretic reformulation: $\G_K(\R)$
contains a subgroup $\C^\times$ 
which acts on $H^{p,q}$ by $z^p \overline{z}^q$.  
If either $K$ or $E$ is in $\R$, as will generally
be the case for us, then 
Hodge symmetry continues to hold.  

If a motive $M$ has Hodge numbers
$\underline{h}^{p,q}$ then the Hodge numbers 
 of its Tate twist $M(j)$ are 
 ${h}^{p-j,q-j} = \underline{h}^{p,q}$.
 The Hodge-normalization of 
 a pure weight motive is the Tate twist
 for which all the nonzero Hodge 
 numbers are in the vector 
 $h=(h^{w,0},\dots,h^{0,w})$ and at least one of
 the outer ones is nonzero.  
  
\para{Zigzag procedure.} 
The procedure we are about to describe is equivalent to a formula
conjectured by Corti and Golyshev \cite{CG} and proved
by different methods in Fedorov~\cite{Fed} and \cite{RV-Mixed}.  The 
procedure is completely combinatorial and only depends on the
interlacing pattern of the roots of $\faa_\infty$ and $\faa_0$ in the
unit circle.

\begin{figure}[htb]
\begin{center}
\begin{tikzpicture}[scale=0.65]
\draw[step=1.0,black, very thin] (0,0) grid (10,3);
\draw[red] (1,3)--(1,3) node[above]{$\frac{1}{4}$};
\draw[red] (2,3)--(2,3) node[above]{$\frac{1}{3}$};
\draw[red] (4,3)--(4,3) node[above]{$\frac{1}{2}$};
\draw[red] (6,3)--(6,3) node[above]{$\frac{2}{3}$};
\draw[red] (7,3)--(7,3) node[above]{$\frac{3}{4}$};
\draw[blue] (0,0)--(0,0) node[below]{$\frac{1}{5}$};
\draw[blue] (3,0)--(3,0) node[below]{$\frac{2}{5}$};
\draw[blue] (5,0)--(5,0) node[below]{$\frac{3}{5}$};
\draw[blue] (8,0)--(8,0) node[below]{$\frac{4}{5}$};
\draw[blue] (9,0)--(9,0) node[below]{$1$};
\draw[red] (10,0)--(10,0) node[right]{$1$};
\draw[red] (10,1)--(10,1) node[right]{$3$};
\draw[red] (10,2)--(10,2) node[right]{$1$};
\draw[blue] (10.5,1)--(10.5,1) node[right]{$1$};
\draw[blue] (10.5,2)--(10.5,2) node[right]{$3$};
\draw[blue] (10.5,3)--(10.5,3) node[right]{$1$};
\draw[red] (-.8,3)--(-.8,3) node[above]{$\alpha_j$'s:};
\draw[blue] (-.8,-.1)--(-.8,-.1) node[below]{$\beta_j$'s:};
\draw[very thick] (0,1) -- (1,0);
\draw[very thick] (1,0) -- (2,1);
\draw[very thick] (2,1) -- (3,2);
\draw[very thick] (3,2) -- (4,1);
\draw[very thick] (4,1) -- (5,2);
\draw[very thick] (5,2) -- (6,1);
\draw[very thick] (6,1) -- (7,2);
\draw[very thick] (7,2) -- (8,3);
\draw[very thick] (8,3) -- (9,2);
\draw[very thick] (9,2) -- (10,1);
\draw[fill=red] (1,0) circle (.1);
\draw[fill=red] (2,1) circle (.1);
\draw[fill=red] (4,1) circle (.1);
\draw[fill=red] (6,1) circle (.1);
\draw[fill=red] (7,2) circle (.1);
\draw[fill=blue] (0,1) circle (.1);
\draw[fill=blue] (3,2) circle (.1);
\draw[fill=blue] (5,2) circle (.1);
\draw[fill=blue] (8,3) circle (.1);
\draw[fill=blue] (9,2) circle (.1);
\end{tikzpicture}
\end{center}
\caption{\label{zigzagfig} The zigzag procedure with input
the family parameter $\fa = {\Phi_2 \Phi_3 \Phi_4}/{\Phi_1 \Phi_5}$
and output the Hodge vector $h=(1,3,1)$}
\end{figure}
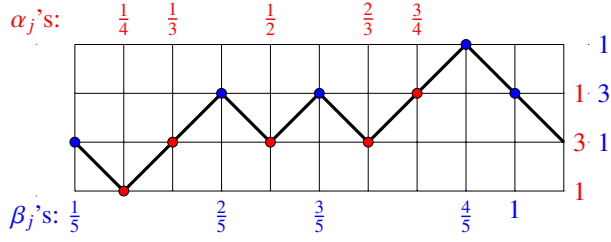

To pass from a family parameter $\fa=\faa_\infty/\faa_0$ to its Hodge vector
$h$ one proceeds as illustrated by Figure~\ref{zigzagfig}.  One
orders the parameters $\alpha_j$ and $\beta_j$, viewed as 
elements in say $(0,1]$; for more immediate readability, we
associate the colors red and blue to $\infty$ and $0$ 
respectively.  One draws a point at $(0,0)$ corresponding to the
smallest parameter in a Cartesian  
plane.   One then proceeds in uniform steps from left to right, 
drawing a point for each parameter and then moving diagonally 
upwards after red points and diagonally downwards after blue points.  
One focuses on one color or the other, counting the number of points 
on horizontal lines.   The numbers obtained form the Hodge  vector
$h$. The red and blue dots yield the same Hodge vector but
  contain   more information. They may be used to describe the
  limiting mixed   Hodge structure at $t=\infty$ and $0$
  respectively. 
  
\para{The completely intertwined case.}  Complete intertwining
of the $\alpha_j$ and $\beta_j$ gives the extreme where the
resulting Hodge vector is just $(n)$.  Beukers and Heckman
\cite{BH} proved that complete intertwining is exactly the condition 
one needs for the monodromy group $\langle h_\infty,h_0 \rangle$
to be finite.  They also established the complete list of such
pairs $(\alpha,\beta)$.  Actually they, like Schwarz 
who famously treated the $n=2$ case more than a century earlier,
worked without our standing assumption $E=\Q$.  Then one 
needs to require complete intertwining of 
all the natural conjugates of $(\alpha,\beta)$ and the list obtained is longer.  

In our setting of $E=\Q$, the corresponding $\gamma$-vectors are of odd lengths~$3$
to~$9$.  There are infinite collections of length $3$ and $5$ given
by coprime positive integers $a,b$: 
\begin{align}
(i)&\quad[-(a+b),a,b], \label{123}  \\
(ii)&\quad[-2(a+b),-a,2a,b,a+b], \nonumber \\
(iii)&\quad[-2a,-2b,a,b,a+b] \nonumber.
\end{align}
Here and always when discussing classification, we omit consideration
of $-\gamma$ whenever $\gamma$ is listed.  
In case $(i)$, the canonical variety consists of just $a+b$ points. 
Removing a variable, the BCM presentation takes the form 
$$
X_{a,b,t}: y^a(1-y)^b-\frac{a^ab^b}{(a+b)^{a+b}}t=0.
$$
For $b=1$ this presentation is already trinomial; in
general, one has to make a non-trivial change of variables
to pass to the trinomial presentation of $X_{a,b,t}$ given by a toric model.     

Beyond the closely related collections $(i)$ - $(iii)$, there 
are only finitely many further $\gamma$, all related
to Weyl groups.  \cite[Table~8.3]{BH} says that,
modulo the quadratic twisting operation $Q(T) \mapsto Q(-T)$, there are  just
one, five, five, and fifteen respectively for the groups
$W(F_4)$, $W(E_6)$, $W(E_7)$, and $W(E_8)$.  
One of the $W(E_6)$ cases is discussed in Section~\ref{dimension}
and the remaining $W(E_n)$ cases are 
similarly treated in \cite{Rob-polys}.  

\para{The completely separated case.} Complete
separation of the $\alpha_j$ and $\beta_j$ gives the
extreme where the resulting Hodge vector is
$(1,1,\dots,1,1)$.  The subcase where
$\faa_0=(T-1)^n$ has the simplifying feature 
that $h_0$ consists of a single Jordan block. 
Families in this subcase have received 
particular attention in the literature; the condition
is sometimes verbalized as MUM, for {\em maximal unipotent monodromy}.    

Classification of families in the completely separated case is 
easier than in the completely intertwined case.  It becomes
trivial in the MUM subcase because $\faa_\infty$ is arbitrary 
except for the fact that it contains no factors of $(T-1)$.  
Accordingly, the number $c_n$ of rank $n$ families
in the MUM subcase is given by a generating function
\begin{eqnarray}
\label{MUMcount}
\sum_{n=0}^\infty c_n x^n & = & \prod_{k=2}^\infty \frac{1}{1-x^{\phi(k)}}  \\
\nonumber &=& 1+x+4 x^2+4 x^3+14 x^4+14 x^5+ \cdots
\end{eqnarray}
Always $c_{2j}=c_{2j+1}$ as under the MUM restriction multiplying
by $(T+1)/(T-1)$ gives a bijection on parameters.  Arithmetic
information about the list underlying $c_4 = 14$ is in \cite{RV}.
For general $n$, the case $\faa_\infty=1+T+\cdots+T^n$ 
is the ``mirror dual'' of the Dwork case discussed after
\eqref{fermat}, and so motives of this family
have been given special attention in the physics literature.   

\para{Signature and the {\em Magma} implementation.}  
A motive defined over a subfield of $\R$ has a signature $\sigma$, which
is the trace of complex conjugation.  For odd weight motives, it
is always zero.  For even weight HGMs $H(\fa,t)$, it depends only on $\fa$ and
the interval $(-\infty,0)$, $(0,1)$, or $(1,\infty)$ in which $t$ lies.  
{\em Magma}'s command \verb@HodgeStructure@ returns 
both the Hodge vector and the signature in coded form.  
To see just the Hodge vector clearly, one can
implement \verb@Q@ as in \eqref{defineq} and 
get the Hodge vector from say
\smallskip

\verb@HodgeVector(HodgeStructure(Q,2));@
\smallskip

\noindent For example, from the gamma vector $[-21,1,2,3,4,5,6]$
one gets the Hodge vector $(1,2,12,2,1)$.

\section{Projective Hypersurfaces}
\label{proj-hypersurf} Here we realize some HGMs in 
the cohomology 
of the most classical varieties of all, 
smooth 
hypersurfaces in projective space.  

\para{Hodge numbers.} 
Let $X \subset \bbP^{\kappa+1}$ be a smooth hypersurface of 
degree $\delta$.  Let $PH^{\kappa}(X,\Q)$ be 
the primitive part of its middle cohomology, meaning the 
part that does not come from the ambient projective space.  If 
$\kappa$ is odd, then this primitive part is all of 
$H^{\kappa}(X,\Q)$. 
If $\kappa$ is even, then the complementary
piece that we are discarding is $\Q(-\kappa/2)$.

Hirzebruch gave a formula for the Hodge numbers
of $PH^{\kappa}(X,\Q)$ as a function
of $\kappa$ and $\delta$.  For example,
 the sum
of the Hodge numbers and first Hodge number
are respectively 
\begin{align}
\label{bettipluri}
b_\kappa & = \frac{(\delta-1)^{\kappa+2}+(-1)^\kappa (\delta-1)}{\delta}, & 
h^{\kappa,0} & =  \binom{\delta-1}{\kappa+1}.
\end{align}
These special cases and Hodge symmetry are sufficient to give
Hodge vectors when $\kappa \leq 3$:
\begin{equation}
\label{projhodge}
{\renewcommand{\arraycolsep}{3pt}
\begin{array}{c|ccc}
\delta  &   \mbox{Curves in $\bbP^2$}  &  \mbox{Surfaces in $\bbP^3$} & \mbox{Threefolds in $\bbP^4$} \\
\hline
    3          &  (1,1) &  (0,6,0) &  (0,5,5,0)   \\
    4          &  (3,3) &  (1,19,1) & (0,30,30,0)   \\
    5          &  (6,6) &  (4,44,4) & (1,101,101,1)  \\
    6          &  (10,10) & (10,85,10) & \, (5, 255, 255,5).  
\end{array}
}
\end{equation}
For $\kappa=1$, either part of \eqref{bettipluri} reduces to the genus formula
for smooth plane curves, $g = (\delta-1)(\delta-2)/2$.

\para{One example for every $(\delta,\kappa)$.}
Let $\delta=e+1\geq 3$ be a desired degree and let 
$\kappa$ be a desired dimension.  Define 
\begin{equation}
\label{gammaekappa}
\gamma=[1,-e,e^2,\ldots,(-e)^{\kappa},(-e)^{\kappa+1}-1,\frac{(-e)^{\kappa+2}+e}{e+1}].
\end{equation}
The toric procedure illustrated by Table~\ref{toroidalchart} yields
the completed canonical variety:
\begin{eqnarray}
\label{smoothhyper} 
X_t:\ \quad u x_{\kappa+1} x_1^e + \sum_{i=2}^{\kappa+2}
  x_{i-1} x_i^e + 
 x_{\kappa+2}^{e+1}  =  0.
\end{eqnarray}
The necessary orthogonality relations on each variable's exponents are illustrated by the 
case of cubic fourfolds where $\gamma=[1,-2,4,-8,16,-33,22]$.  Then  
\eqref{smoothhyper} becomes 
\[
X_t:  u x_5 x_1^2 + x_1x_2^2+x_2x_3^2+x_3x_4^2+x_4x_5^2+x_5x_6^2+x_6^3=0.
\]
For $x_5$, the relation is that $m_{5*}=(1,0,0,0,2,1,0)$ is
orthogonal to $\gamma$.   In general, partial derivatives of \eqref{smoothhyper} are very simple since the row vectors
$m_{i*}$ have just two nonzero entries, except for the case $i=\kappa+1$ 
and its three nonzero entries.   It is then a pleasant exercise to check via the Jacobian criterion that $X_t$ is 
smooth for $t \in \C^\times-\{1\}$. 

The degree of the rational function $\fa$ determined by $\gamma$ can be computed 
uniformly in $(\delta,\kappa)$ as the cancellations to be analyzed 
are very structured.  This degree agrees with the Betti 
number $b_\kappa$ from \eqref{bettipluri}.   Thus
$H(\gamma,t)$ is the full primitive middle cohomology
of $X_t$, while {\em a priori} it might have been a proper
subspace.   The zigzag procedure for computing 
Hodge numbers must agree in the end with the Hirzebruch
formula.  The reader might want to check the above case of 
cubic fourfolds, where Hirzebruch's full formula gives $(0,1,20,1,0)$.

\para{All examples for a given $(\delta,\kappa)$.}  An interesting
problem is to find all $\gamma$ which give projective smooth
$\kappa$-folds of degree $\delta$.  For small parameters, this
problem can be solved by direct computation.  For example,
consider $(\delta,\kappa)=(3,4)$, thus cubic fourfolds.  
In this case, one has the standardization $[-33,-8,-2,1,4,16,22]$ of
the above example, and then exactly ten more:  
\[
{\renewcommand{\arraycolsep}{0pt}
{
\!\! \begin{array}{ll}
\, [-48, -15, -12, 5, 16, 24, 30], \; \; \; & [-36, -9, -4, 3, 8, 18, 20],\\
\, [-48, -12, -3, 1, 6, 24, 32], & [-33, -16, -4, 2, 8, 11,32],\\
\, [-48,-12, -3, 6, 16, 17, 24], & [-33, -10, -7, 5, 11, 14, 20],\\
\, [-36, -16, -9, 3, 8, 18, 32], & [-33,-4, -1, 2, 8, 11, 17], \\
 \,[-36, -9, -8,4, 15, 16, 18], &  [-21, -20, -16, 7, 8, 10, 32].
 \end{array}
\!\! }
}
\]

\section{Dimension reduction}
\label{dimension}
    An HGM $H(\gamma,t)$ is defined in terms of a $\kappa$-dimensional variety
but its Hodge vector $(h^{w,0},\dots,h^{0,w})$ raises the question  
of whether it also comes from a variety of dimension $w = \kappa-2j$.  
The exterior zeros for low degree projective hypersurfaces  
as illustrated in \eqref{projhodge} raise the same question.  The
generalized Hodge conjecture says that this dimension 
reduction is always possible.  We illustrate here some 
of the appealing geometry that arises from 
reducing dimension.

\para{Reduction to points.} When $w=0$ the reduction 
to dimension zero is possible in all cases.  For example, 
$\gamma = [-12,-3,1,6,8]$ corresponds to entry 45
on the Beukers-Heckman list \cite[Table~8.3]{BH}.
Formula \eqref{toroidal} then gives a family $X_t$ of cubic surfaces.
An equation whose roots correspond to the famous twenty-seven lines
on $X_t$ is 
\begin{align}
\label{xt}
2^4tx^3(x^2-3)^{12} - 3^9(x^3-3 x^2+x+1)^8(x-2) = 0.
\end{align}
The Galois group of this polynomial $g(t,x)$ 
for generic $t \in \Q^\times  - \{1\}$ 
is $W(E_6)$.  It has 
$51840 = 2^7 3^4 5$ elements and is also the monodromy group $\Gamma = \langle h_\infty,h_0 \rangle$.

\para{Reduction via splicing.}  Suppose $\gamma$ can be written as the concatenation
of two lists each summing to zero.  Then one can use a general splicing technique from \cite[\S6]{BCM} to
reduce the dimension by two.  This technique is behind the scenes even of our introduction: 
in the family of examples there, the canonical varieties for $[-2,-2,1,1,1,1]$ are three-dimensional, 
although the more familiar source varieties are just the Legendre curves \eqref{legendre-equation}.

For an example complicated enough to be representative of the general case,
take $\gamma = [-12,-3,-2,1,1,1,6,8]$ 
so that the canonical variety has dimension $\kappa = 5$.  
The Hodge vector is just $(3,3)$, so one would like to 
realize $H(\gamma,t)$ in the cohomology of a curve.  

Splicing is possible because both $[-12,-3,1,6,8]$ and  $[-2,1,1]$ sum to zero.  No further
splicing is possible, but fortunately we have just treated the first sublist
 by other means.  
Splicing corresponds to taking a fiber product over the $t$-line
which in turn corresponds to just multiplying rational functions. 
In our case, solving \eqref{xt} for $t$ to get the first factor,
the dimension-reduced variety is given by
 \begin{equation}
 \label{spliceprod}
 \frac{3^9(x-2) \left(x^3-3 x^2+x+1\right)^8}{2^4 x^3
    \left(x^2-3\right)^{12}} \cdot \frac{2y-1}{y^2} = t.
 \end{equation}
 The variable $y$ from $[-2,1,1]$ enters only quadratically and so \eqref{spliceprod} defines
 a double cover of the $x$-line.   Taking the discriminant with respect 
 to $y$ and removing unneeded square factors presents this hyperelliptic
 curve in standard form:
\[
z^2 = -3 (x-2) g(t,x).
\]
As the right side has degree 28, this curve has 
genus $13$.   
  
    In both the new examples of this section, the middle cohomology
of the dimension-reduced varieties contains not only the
desired motives, with Hodge vectors $(6)$ and $(3,3)$ respectively,
but also parasitical motives, with Hodge vectors $(21)$ 
and $(10,10)$.    In this regard, they are less attractive
than the original canonical varieties.  HGMs
provide many illustrations like these two of the motivic principle that
a motive $M$ comes from many varieties $X$, 
and often no single $X$ should be viewed as the best
source.  

\section{Distribution of Hodge vectors}
\label{distribution}
In this section, we explain one of the great features of HGMs: they
represent many Hodge vectors.  

\para{Completeness in ranks $\leq 19$.}
By direct computation starting from all family
parameters $\fa$ in degrees $\leq 19$, we have
verified the following fact.   
{\it Let $h=(h^{w,0},\cdots,h^{0,w})$ be a vector of positive integers
  satisfying $h^{q,p}=h^{p,q}$ for all $p+q=w$ and let
  $n=\sum_{i=0}^w h^{p,w-p}$. Then if $n\leq 19$ there exists an
  HGM with Hodge vector $h$.}
  
\para{Many families per Hodge vector in ranks $\leq 100$.} 
In ranks $20$ to $23$, the only vectors not realized by a family
of HGMs are
\[
\begin{array}{rc}
20: & (6,1,1,1,2,1,1,1,6), \\
22: & (6,1,1,1,1,2,1,1,1,1,6), \\
22: & (4,1,2,1,1,1,2,1,1,1,2,1,4), \\
23: & (1,21,1).
\end{array}
\]
Table~\ref{rank24} gives a fuller sense of the situation for
$n=24$, where there are about $460,000,000$ family parameters.
It gives the 
extremes of the list of $4096$ possible Hodge
vectors $h$, sorted by how many families realize $h$. 
\begin{table}[htb]
\[
\begin{array}{c|r}
h & \# \\
\hline
$(9, 1, 1, 2, 1, 1, 9)$ & 0 \\
$(7, 1, 1, 1, 1, 2, 1, 1, 1, 1, 7)$ & 0 \\
$(1, 6, 1, 1, 1, 1, 2, 1, 1, 1, 1, 6, 1)$& 0 \\
$(4, 1, 3, 1, 1, 1, 2, 1, 1, 1, 3, 1, 4)$& 0 \\
$(5, 1, 2, 1, 1, 1, 2, 1, 1, 1, 2, 1, 5)$& 0 \\
$(6, 1, 1, 1, 1, 1, 2, 1, 1, 1, 1, 1, 6)$& 0 \\
$(4, 1, 1, 2, 1, 1, 1, 2, 1, 1, 1, 2, 1, 1, 4)$& 0 \\
$(4, 1, 2, 1, 1, 1, 1, 2, 1, 1, 1, 1, 2, 1, 4)$& 0 \\
$(6, 2, 1, 1, 1, 2, 1, 1, 1, 2, 6)$&   $2$\\
$(8, 1, 1, 1, 2, 1, 1, 1, 8)$&  $4$\\
$(1, 22, 1)$&$	                           4$\\
$(8, 1, 1, 4, 1, 1, 8)$&$	                   6$\\
$\vdots$&$\vdots$\\
$(1, 5, 6, 6, 5, 1)$&$	              7637828$\\
$(1, 2, 4, 5, 5, 4, 2, 1)$&$	      7982874$\\
$(2, 4, 6, 6, 4, 2)$&$	              9504072$\\
$(1, 4, 7, 7, 4, 1)$&$	              9905208$
\end{array}
\]
\caption{\label{rank24} Hodge vectors with total $24$ and
their number of hypergeometric realizations}
\end{table}

The ratio of the numbers just reported say 
that the number of family parameters per Hodge
vector in degree $24$ is about $113,000$.  
This ratio increases to a maximum at
$n=58$ where it is about four million.   
It then decreases to zero, with some 
approximate sample values being
two million for $n=100$ but only 
$0.00001$ for $n=300$.   These 
numbers are computed via generating 
functions, similar to \eqref{MUMcount}
but more complicated.

 \para{Perspective.}
 Section~\ref{proj-hypersurf} offers some perspective on the 
 general inverse problem of finding an irreducible 
 motive $M \in \cM(\Q,\Q)$ with a given Hodge
 vector.  
 From \eqref{bettipluri}-\eqref{projhodge}, one sees that the 
 Hodge vectors coming from hypersurfaces are very sparse.
 When one looks at broader standard classes of varieties,
 such as complete intersections in projective spaces,
 more Hodge vectors arise, but they all have the same 
 rough form: bunched in the middle.  {\em Ad hoc} 
 techniques, such as reducing Hodge numbers
 by imposing singularities, give many more Hodge vectors. 
 But for many $h$, it does not seem easy to find
 a corresponding motive and then prove irreducibility in this geometric way. 
  For example, imposing $k$ ordinary double 
 points on a sextic surface reduces the
 Hodge vector to $(10,85-k,10)$.  However
 the family of sextic surfaces is only $68$-dimensional,
 and so it would it seem to be difficult to get down to 
 e.g.\ $(10,1,10)$.  There does not seem
 to be even a conjectural expectation of 
 which Hodge vectors arise
 from irreducible motives in 
 $\cM(\Q,\Q)$.

\para{The cases $(1,b,1)$.}  
One could go into much more detail about 
the families behind any given Hodge vector. 
Here we say a little more about the
cases $(1,b,1)$, which 
are particularly interesting for
several reasons.  The $\gamma$ giving Hodge vectors
of the form $(1,b,1)$  typically have canonical dimension 
$\kappa = \dim(X_{\gamma,t})$ greater than 
two, posing instances of the dimension reduction problem.  
If $b \leq 19$, then the moduli theory
of $K3$ surfaces says that there is at least one
family $Y_{t}$ of $K3$ surfaces
also realizing $H(\gamma,t)$.  Finding
such a family is a challenge.

Cases with $b \geq 20$ present a greater 
challenge, as they cannot be realized by K3 surfaces. 
There are seventy-two parameters giving 
$(1,20,1)$.  None of the eleven 
listed in Section~\ref{proj-hypersurf} 
can be spliced, underscoring the
difficulty of dimension reduction.   
One of the four gamma vectors 
giving $(1,22,1)$ has canonical dimension
eight, namely $[-60, -5, -4, -3, -2, 8, 9, 10, 12, 15, 20]$.
The other three have canonical dimension
ten:
$$
\begin{array}{l}
\, [-66, -11, -6, -5, -4, -4, 1, 2, 8, 12, 18, 22, 33],\\
\,  [-60, -15, -9, -6, -4, -2, 3, 5, 8, 12, 18, 20, 30],\\
\, [-33, -10, -6, -4, -4, -1, 2, 2, 5, \; \, 8, 11, 12, 18].
\end{array}
$$
In all four cases, there are many ways to splice, but
no path to a surface.

\section{Special and semi HGMs}
\label{special-motive}
We have so far been excluding the singular specialization point $t=1$ from consideration.
  Now we explain how it yields a particularly interesting 
  motive $H(\fa,1) \in \cM(\Q,\Q)$. 
We also explain how other interesting motives arise when the
family parameter $\fa$ is reflexive, in the sense of satisfying $\fa(-T)=\fa(T)^{-1}$.

\para{Interior zeros.} A Hodge-normalized motive $M \! \in \! \cM(\Q,\Q)$ of weight $w$
has Hodge vector $h = (h^{w,0},\dots,h^{0,w})$ with $h^{w,0}=h^{0,w}>0$.   
But for the Hodge vectors explicitly considered so far, the remaining
numbers $h^{p,w-p}$ are also positive.   There is a reason for this restriction: 
Griffiths transversality says that any collection of motives moving
in a family with irreducible monodromy group has Hodge vector with
no interior zeros.   Special and semi HGMs do not move in families,
and they include cases with interior zeros.  

\para{Special HGMs.}  The way to account for the double point
on the canonical variety $X_{\gamma,1}$ is to first of all take
inertial invariants with respect to the monodromy operator 
$h_1$.  In the orthogonal case, this already give
the right motive $H(\fa,1)$.   Its Hodge vector
differs from the generic Hodge vector only in that
$h^{w/2,w/2}$ is decreased by $1$.   In the 
symplectic case, the motive of inertial invariants
is mixed, and quotienting out by its submotive
of weight $w-1$ and rank $1$ gives $H(\fa,1)$.  
Its Hodge vector now comes from the generic one
by decreasing the two central Hodge numbers by
$1$.   These drops obviously can cause interior
zeros, as in $(10,1,10) \rightarrow (10,0,10)$ 
or $(1,1,1,1,1,1) \rightarrow (1,1,0,0,1,1)$.

\para{Semi HGMs.}  
For a reflexive parameter $\fa$ and any $t \in \Q^\times$, 
the motives $H(\fa,t)$ and $H(\fa,t^{-1})$
are quadratic twists of one another.  
The interest in reflexive parameters is that 
nongeneric behavior is thereby forced at $t=\pm 1$.   
The motive  
$H(\fa,(-1)^n)$ is a direct sum of 
two motives  in $\cM(\Q,\Q)$ of roughly equal
rank.  We call the summands semi HGMs and their Hodge vectors
can have many interior zeros.   For example, the summands of 
$H(\Phi_2^{16}/\Phi_1^{16},1)$ 
are studied in \cite{Rob} and the two Hodge vectors
are 
\begin{equation}
\label{dyadic16}
\begin{array}{c}
\, (1,0,1,0,1,0,1,0,0,1,0,1,0,1,0,1), \\
\, (1,0,1,0,1,0,0,0,0,1,0,1,0,1).
\end{array}
\end{equation}
There is a similar 
decomposition of $H(\fa,-(-1)^n)$, but only
after viewing it in $\cM(\Q(i),\Q(i))$.

\section{Point counts}
\label{sect:arithmetic}

We now turn to arithmetic.  The point counts of this section form the
principal raw material from which the $L$-functions studied 
in the remaining sections are built.

\para{Background.} 
Let $X$ be a smooth projective variety
over $\Q$.  Then for all primes $p$ outside a finite set $S$, the
equations defining $X$ have good reduction and so define
a smooth projective variety over $\F_p$.  For any power $q=p^e$, 
one has the finite set of solutions $X(\F_q)$ to the defining
equations.  The key invariants that need to be input into 
the motivic formalism are the cardinalities $|X(\F_q)|$,
and famous results of Grothendieck, Deligne, and 
others provide the tools.  

The vector spaces $H^k(X,\Q)$ do not see that $X$ is 
defined over $\Q$.  The arithmetic origin of $X$ 
yields extra structure as follows.  For any prime $\ell$, one
can extend coefficients to obtain vector spaces 
$H^k(X,\Q_\ell)$ over the field $\Q_\ell$ of $\ell$-adic
numbers.  Then the group $\Gal(\overline{\Q}/\Q)$ 
acts on $H^k(X,\Q_\ell)$.   
 
For every prime $p$ the group $\Gal(\overline{\Q}/\Q)$ 
contains Frobenius elements $\Fr_p$, well-defined up to 
ambiguities that will disappear from our considerations.  
For any power $q=p^e$ of a prime $p \not \in S$, and any $\ell \neq p$, one
has the trace of the operator $\Fr_{q} = \Fr_p^e$ acting on 
$H^k(X,\Q_\ell)$.  These $\ell$-adic 
numbers are in fact rational and independent of 
$\ell$.   We emphasize the independence of $\ell$ by denoting them $\Tr(\Fr_q|H^k(X,\Q))$.  
The connection with point counts is the Lefschetz trace formula:
$|X(\F_q)| = \sum_k (-1)^k \Tr(\Fr_q|H^k(X,\Q))$.  The
left side for fixed $p$ and varying $e$ determines the summands
on the right side in principle because the complex eigenvalues of
$\Fr_p$ on weight $k$ cohomology have absolute value 
$p^{k/2}$.

Much of this transfers formally to the motivic setting. 
Thus for a motive $M$ and a prime $\ell$, there is
an action of $\Gal(\overline{\Q}/\Q)$ on the corresponding
$\ell$-adic vector space $M_\ell$.   This action has image in 
$G_M(\Q_\ell)$.  Indeed the Tate conjecture predicts 
that the $\Q_\ell$-Zariski closure of the image of $\Gal(\overline{\Q}/\Q)$ is all of
$G_M(\Q_\ell)$. 

One technical problem with Andr\'e's category 
$\cM(\Q,\Q)$ is that the projectors used to define 
motives are not known to come from algebraic cycles.  As a consequence,
for a general $M \in \cM(\Q,\Q)$ the above compatibility of Frobenius traces
is not known.  However this problem does not arise for 
hypergeometric motives, because they are essentially
the entire middle cohomology of varieties.  Accordingly
one has well-defined rational numbers $\Tr(\Fr_q|H(\gamma,t))$. 
There are similar technical problems at the 
primes $p \in S$, but they do not affect our
computations and we will ignore them.

\para{Wild, tame, and good primes.}  Returning now to
very concrete considerations, we sort primes for a 
parameter $(\gamma,t)$ as follows.  
  A prime $p$
is {\em wild} if 
it divides a $\gamma_j$.
For $t \neq 1$, a prime $p$ is {\em tame} if 
 it
is not wild but it divides either the numerator of $t$, the denominator
of $t$, or the numerator of $t-1$; these last three conditions say that $t$ is $p$-adically close to the
special points $0$, $1$, and $\infty$ respectively.   For $t=1$,
no primes are tame.   We say that a prime is {\em bad}
if it is either wild or tame, and all other primes are {\em good}.

\para{Split powers of a good prime $p$.}   
A power $q$ of a good prime $p$ is {\em split} for $\gamma$ if $q\equiv 1 \bmod m$,
where $m$ is the least common multiple of the $\gamma_j$.   One then 
 has a collection of Jacobi sums indexed by characters $\chi$ of $\F^\times_q$:
\[
J(\gamma,\chi) := \prod_{j=1}^n g(\omega^{\alpha_j (q-1)} \chi,\psi) \overline{g(\omega^{\beta_j (q-1)} \chi,\psi)}.
\]
Here $\psi: \F_q \rightarrow \C^\times$ is any nonzero additive character, $\omega: \F_q^\times \rightarrow \C^\times$ is
any generator of the group of multiplicative characters, $(\alpha,\beta)$ underlies $\gamma$ as in
Section~\ref{monodromy}, and $g(\rho,\psi)=\sum_{t\in\F_q^\times} \rho(t)\psi(t)$ is the standard Gauss sum.
The desired quantity is then given by a
sum 
due to Katz \cite[p.\ 258]{Katz-ESDE}.  Renormalizing to fit our conventions, it
is 
\begin{equation}
\label{hyperg-trace}
\Tr(\Fr_q|H(\gamma,t)) = 
\frac {q^{\phi_0}}{1-q}\sum_\chi
\frac{J(\gamma,\chi)}{J(\gamma,1)} 
\,\chi(t).
\end{equation}
Here $\phi_0$ is the vertical coordinate
of a lowest point on the zigzag diagram of $\gamma$, e.g.\ $\phi_0=-1$ in
Figure~\ref{zigzagfig}.

\para{General powers of a good prime $p$.} 
The Gross-Koblitz formula lets one 
replace the above Gauss sums by values
of the $p$-adic gamma function.  This is both 
a computational improvement and extends the formula to 
all powers of any good prime.  
With this method, the desired integers $\Tr(\Fr_q|H(\gamma,t))$ 
are first approximated $p$-adically.  Errors are under control and
exact values are determined from sufficiently good approximations. 
See \cite{BCM} 
for a closely related approach to the essential numbers $\Tr(\Fr_q|H(\gamma,t))$ 
and references to earlier contributions.   

\section{Frobenius polynomials} 
\label{frobenius}
 Frobenius 
polynomials are a concise way of packaging
the point counts of the preceding section.  
They play the leading role in the formula
for $L$-functions of the next section.   
After saying what they are, this section explains several reasons 
why they are useful, even before 
one gets to $L$-functions.

\para{Capturing point counts.}
Consider the numbers $c_{p,e} = \Tr(\Fr^e_p|M) \in \Q$
 for a fixed motive $M \in \cM(\Q,\Q)$ of rank $n$, a fixed good prime $p$, and varying $e$.  
They can be captured in a single degree $n$ polynomial
$F_p(M,x)=\det(1-\Fr_p x|M)$.  The relation, which comes from summing
the geometric series belonging to each of the $n$ eigenvalues, is
\begin{equation}
\label{expdef}
\exp \left(\sum_{e=1}^\infty \frac{c_{p,e}}{e} x^e \right) = \frac{1}{F_p(M,x)}.
\end{equation}
Write 
\[
F_p(M,x)=1 + a_{p,1} x + \dots + a_{p,n-1} x^{n-1} + a_{p,n} x^n.
\]
Then the $c_{p,e}$ for $e \leq k$ determine $a_{p,k}$. 
Thus the $c_{p,e}$ for $e \leq n$ determine $F_p(M,x)$.
But, even better, Poincar\'e
duality on a source variety ultimately implies that
one has $a_{p,e} = \epsilon(p) a_{p,n-e} p^{(n-2e)w/2}$
for a sign $\epsilon(p)$.   For HGMs, this sign is known
and in fact always $1$ when $w$ is odd.  So 
$F_p(M,x)$ can be computed using only $c_{p,e}$ for 
$e \leq \lfloor n/2 \rfloor$.

\para{Relation with Hodge vectors.}  Indexing by weight $w$, consider as examples the 
rank six family parameters
\begin{align}
\fa_0 & = \frac{\Phi_3 \Phi_{12}}{\Phi_1 \Phi_2 \Phi_8}, &
\fa_1 & = \frac{\Phi_3 \Phi_{12}}{\Phi_1^2 \Phi_8}, &
\fa_5 & = \frac{\Phi^3_3}{\Phi_1^6}.  \label{qw}
\end{align}
The first two are the families from Section~\ref{dimension}, with
Hodge vectors respectively $(6)$ and $(3,3)$; the last one has
Hodge vector $(1,1,1,1,1,1)$.  Specializing at a 
randomly chosen common point gives motives
$M_{6,w} = H(\fa_{w},3/2)$.

After the required initialization of a variable $x$ by \verb@_<x>:=PolynomialRing(Integers())@,
and after inputting \verb@Q@$w$ as in \eqref{defineq},
{\em Magma} quickly gives some Frobenius polynomials
via e.g.\ \verb@EulerFactor(Q0,3/2,5)@:
\begin{align*}
F_5(M_{6,0},x) & =  1-x-x^5+x^6, \\
F_7(M_{6,0},x) & =   1 \w\w\w\w\w\w\w-x^6, \\
& \\
F_5(M_{6,1},x) & = 1+\w x + \w 6 x^2 + 16 x^3+ \cdots, \\
F_7(M_{6,1},x) & =  1-2x + 12 x^2 - 28 x^3 + \cdots, \\
&\\
F_5(M_{6,5},x) & =  1-\w 9 x + 5 \ccdot 156 x^2 - 5^3 \ccdot 2556 x^3  + \cdots, \\
F_7(M_{6,5},x) & =  1+12 x + 7 \ccdot 888 x^2 + 7^3 \ccdot 1816  x^3 + \cdots.
\end{align*}
These displays illustrate a basic motivic principle: as weight 
increases, motives of a given rank $n$ become more
complicated.   A more refined principle involves Hodge numbers   
and can be expressed by forming 
a weakly increasing vector $(s_1,\dots,s_n)=(0,\dots,w)$, 
where an entry $i$ appears $h^{i,w-i}$ times.  Then the Newton-over-Hodge
inequality is   
$
\ord_p(a_{p,k}) \geq \sum_{j=1}^k s_j.
$
For $k=1$, \dots, $6$, these lower bounds from the Hodge
vector $(1,1,1,1,1,1)$ controlling $M_{6,5}$ are 
$(0,1,3,6,10,15)$.  For $(3,0,0,0,0,3)$ the 
bounds  $(0,0,0,5,10,15)$ would be smaller, 
leaving more possibilities for Frobenius polynomials.
In this sense, spread out Hodge vectors correspond to
more complicated motives.  

\para{Congruences.} Reduced to $\F_\ell$, the numbers
$a_{p,k}$ for $p \neq \ell$ depend only on the mod $\ell$ Galois 
representation belonging to $M$.  In our 
examples, suppose one kills $\ell = 2$ in \eqref{qw}
by replacing all $\Phi_{2^j a}$ by $\Phi_a^{\phi(2^j)}$.  
Then $\fa_0$ and $\fa_1$ both become $\fa_5$.  
This agreement implies that $F_p(M_{6,w},x) \in \F_\ell[x]$ is
independent of $w$.   This independence  can
be seen for the primes $5$ and $7$  in the displayed
Frobenius polynomials.  The analogous
congruences hold for any $\ell$, when
one changes our Tate twist convention
to make the weight of $H(\fa,t)$ the 
number of integers among the $\alpha_j$ and $\beta_j$, 
minus one.  This web of congruences, like the web corresponding
to splicing considered in Section~\ref{dimension}, makes it clear that 
HGMs constitute a natural collection of motives.

\para{Finite Galois groups.}  Frobenius polynomials render Galois-theoretic 
aspects of the situation very concrete.   As a warm-up, consider $\fa_0$ 
as a representative of the relatively familiar case of 
 ordinary Galois theory.  Here the $\ell$-adic representations all come from
a single representation $\Gal(\overline{\Q}/\Q) {\rightarrow} W(E_6) \subset GL_6(\Q)$.
Let $\lambda_p$ be the partition of $27$ obtained by taking the
degrees of the irreducible factors of $g(3/2,x)$ from \eqref{xt}.  Then
the twenty-five possibilities for the pair $(\lambda_p,F_p)$ 
correspond to the twenty-five conjugacy classes in the
finite group $W(E_6)$.   If one can collect 
enough classes, then one can conclude that 
the image $G$ is all of $W(E_6)$.  In our example $t=3/2$, 
the above primes $5$ and $7$ give $(5^5 1^2,1-x-x^5+x^6)$ and  $(6^4 3,1-x^6)$ respectively.
In ATLAS notation, these are the classes $5A$ and $6I$.  They do 
not quite suffice to prove $G=W(E_6)$.  But the prime $11$ gives the class $12C$ and 
since no maximal subgroup contains elements from $5A$, $6I$, and $12C$, indeed
$G=W(E_6)$.    

The Chebotarev density theorem
says that each pair appears proportionally to the
number of elements in its conjugacy class.  For example, 
the classes $5A$, $6I$ and $12C$ occur with frequency 
$1/10$, $1/12$, and $1/12$ respectively.   

\para{Infinite Galois groups.} 
The cases $\fa_1$ and $\fa_5$ are beyond classical
Galois theory as the motivic Galois groups have 
positive dimension.  But the situation remains quite 
similar.   Consider for example odd weight motives of 
rank $n=2r$ so that $G$ is in the conformal symplectic group
$\CSp_n$.   The Weyl group of $\CSp_n$ is the hyperoctahedral
group $W(C_r)$ of signed permutation matrices, with
order $2^r r!$.    
A separable $F_{p}(M,x)$, being conformally palindromic,
has Galois group within $W(C_r)$.   If it has
Galois group all of $W(C_r)$ then $G$ necessarily 
contains a certain twisted maximal torus.  Suppose a second
prime $p'$ satisfies the same condition and
moreover the joint Galois group of $F_p(M,x) F_{p'}(M,x)$
is all of $W(C_{r}) \times W(C_r)$.  Then $G$ contains
two maximal tori which are sufficiently different to force
$G = \CSp_n$, by the classification of subgroups 
containing a maximal torus.   

To analyze a given motive, the necessary computations can be 
done using {\em Magma}'s \verb@GaloisGroup@ 
command.   The order of the Galois group of $F_{p}(M_{6,1},x)$ 
is $16$, $16$, $4$, $48$, $48$ for $p=5$, $7$, $11$, $13$, 
$17$, and the pair $(p,p') = (13,17)$ satisfies 
the criterion.   For $5 \leq p < 100$, all $F_p(M_{6,5},x)$ have
Galois group $W(C_3)$ except $p=13$.  Excluding 
13, all $\binom{22}{2} = 231$ pairs $(p,p')$ satisfy
the criterion.   In general, it becomes easier to establish
genericity as the weight increases, a reflection of the 
growth in complexity discussed above.

Applying this two-prime technique to the special and semi HGMs 
of Section~\ref{special-motive} suggests that almost always their motivic Galois
groups are as big as possible.  In particular, the exotic Hodge vectors with
interior zeros arising there indeed come from irreducible motives.  
Details in the case \eqref{dyadic16} are given in \cite{Rob}.  

The Chebotarev density theorem extends to the full motivic
setting if all the $L$-functions described below 
have their expected analytic properties.  
Readers wishing to see a glimpse of this theory 
can compute hundreds of $a_{p,1}/p^{w/2}$ 
for $M_{6,w}$ for $w=1$ or $w=5$.  
 By all appearances, the data matches the Sato-Tate
measure $\mu$, meaning the pushforward
of Haar measure on the compact group $\Sp_6$ to $[-6,6]$ via 
the defining character.    One would have to 
compute thousands of $a_{p,1}/p^{w/2}$ before one 
could confidently distinguish this measure
from the Gaussian measure
of mean $0$ and standard deviation $1$.  

\section{$L$-functions}
\label{sect:L-functions}

We now finally define $L$-functions and illustrate
how everything works by some numeric computations.  

\para{Local invariants.} 
  Let $M \in \cM(\Q,\Q)$ be a motive of rank $n$ and weight $w$, having bad
reduction within a finite set $S$ of primes.  We have
discussed two types of local invariants associated to $M$. 
Corresponding to the place $\infty$ of $\Q$ is the
Hodge vector $h = (h^{w,0},\dots,h^{0,w})$ with
total $n$, and also a signature $\sigma$. 
Corresponding to a prime $p \not \in S$ is 
the degree $n$ Frobenius polynomial 
$F_p(M,x)$.   For primes $p \in S$, there
is also a Frobenius polynomial $F_p(M,x)$, now 
of degree $\leq n$, and moreover a conductor exponent $c_p \geq n-\mbox{deg}(F_p(M,x))$,
both to be discussed shortly.  The 
conductor of $M$, which can be 
viewed as quantifying the severity
of its bad reduction, is the integer $N = \prod_{p \in S} p^{c_p}$.  

 \para{Formal products.}  The local invariants can be combined into 
a holomorphic function in the right half-plane
$\Real(s)>\frac{w}{2}+1$, called the completed
$L$-function of $M$:
\begin{equation}
\label{Lprod}
\Lambda(M,s) = N^{s/2} \Gamma_{h,\sigma}(s) \prod_p \frac{1}{F_p(M,p^{-s})}.
\end{equation}
The product over primes alone is the 
$L$-function $L(M,s)$, while the remaining factors give the 
completion.  The infinity factor is given by an explicit 
formula: 
\begin{equation} 
\label{inffact}
\Gamma_{h,\sigma}(s) = \Gamma_\R(s-\frac{w}{2})^{h_+} \Gamma_\R(s-\frac{w}{2}+1)^{h_-} \prod_{p<q} \Gamma_\C(s-p)^{h^{p,q}}  \!\! .
\end{equation}
Here 
$\Gamma_\R(s)  = \pi^{-s/2} \Gamma(s/2)$ and
$\Gamma_\C(s) = 2 (2 \pi)^{-s} \Gamma(s)$.
 The factors involving $h_\pm = (h^{w/2,w/2} \pm (-1)^{w/2} \sigma)/2$ only appear when $w$ is even; in the common case that
$\sigma = 0$, they can be replaced by $\Gamma_\C(s-\frac{w}{2})^{h^{w/2,w/2}/2}$, by the duplication
formula.  

Both the $L$-function and the completing factor are
multiplicative in $M$ so that $\Lambda(M_1 \oplus M_2,s) = \Lambda(M_1,s) \Lambda(M_2,s)$.
Another simple aspect of the formalism is that Tate twists correspond
to shifts: $\Lambda(M(j),s) = \Lambda(M,s+j)$.  

\para{Expected analytic properties.}     The $L$-function $L(\Q,s)$ of the unital motive $\Q$ is just the Riemann zeta function
 $\zeta(s) = \prod_{p} (1-p^{-s})^{-1}$, and the completing factor is $\Gamma_\R(s)$.  
 Riemann established that $\Lambda(\Q,s)$ has a meromorphic continuation to the
 whole $s$-plane, with poles only at $0$ and $1$; moreover he proved that 
 $\Lambda(\Q,1-s) = \Lambda(\Q,s)$. The product $\Lambda(M,s)$ is expected to have similar analytic properties.   First, for $M$
irreducible and not of the form $\Q(j)$, there should be an analytic continuation
to the entire $s$-plane, bounded in vertical strips.  Second, always
\begin{equation}
\label{fctnl-eqn}
\Lambda(M,w+1-s)=\epsilon\Lambda(M,s),
\end{equation}
for some sign $\epsilon$.   For comparison with Section~\ref{automorphy},
note that most everything said in the last three sections generalizes
to motives in $\cM(\Q,E)$, with Frobenius polynomials being in $E[x]$.   However 
\eqref{fctnl-eqn} takes 
the more complicated form $\Lambda(\overline{M},w+1-s)=\epsilon\Lambda(M,s)$,
with $\overline{M}$ the complex conjugate motive and $\epsilon$ only on the unit circle.  

 \para{Determining invariants at bad primes.}   One approach to the conductor exponents $c_p$ and 
 Frobenius polynomials $F_p(M,x)$ associated to 
 bad primes $p$ is to compute them directly by 
 studying the bad reduction of an
 underlying variety.   For an HGM $H(\fa,t)$, {\em Magma} takes this approach for 
 primes which are tame for $(\fa,t)$, as sketched in Section~\ref{bad}.

    A very different approach  
 uses the fact that the list of possible 
 $(c_p,F_p(M,x))$ for a given prime $p$ is finite, 
and the product \eqref{Lprod} has the conjectured
analytic properties for at most one member 
of the product list.  The current state of HGMs 
for the wild primes of $\fa$ mixes
the two approaches: we first greatly reduce the 
length of the lists by using proved and
conjectured general facts.  Then we 
search within the much smaller product list for 
the right quantities. 

   Our view is that numerical computations such as 
those that follow in this section and Section~\ref{numerical}
admit only one plausible interpretation: the bad factors
have been properly identified and the 
analytic properties indeed hold.   However 
rigorous confirmation does not seem to be
in sight at the moment, despite the progress described
in Section~\ref{automorphy}.    
   
\para{A rank four example.} For $\fa$ of degree $\leq 6$ and $t=1$, 
Watkins numerically identified all the bad quantities, 
so that the corresponding $L$-functions
are immediately accessible on {\em Magma}.
For example, take the family parameter 
to be $\fa=\Phi_2^2 \Phi_{12}/\Phi_{18}$, implemented
as always by modifying \eqref{defineq}.
At the specialization point $t=1$, the
Hodge vector is $(1,1,1,1)$.  The corresponding $L$-function, set up
so that calculations are 
done with 10 digits of precision, is
\smallskip

\verb@L := LSeries(Q,1:Precision:=10);@
\smallskip

\noindent The bad information stored in {\em Magma} is revealed by
\verb@EulerFactor(L,@$p$\verb@)@ and \verb@Conductor(L)@ 
to be $F_2(M,x)=1+2x$, $F_3(M,x)=1$,
and $N = 
2^6 3^9$.   The sign $\epsilon$
is calculated numerically, with \verb@Sign(L)@ 
returning $-1.000000000$.  So the order of vanishing 
of $L(M,s)$ at the central point $s=2$ should be odd.  
This order is apparently three since
\smallskip

\verb@Evaluate(L,2:Derivative:=1);@
\smallskip

\noindent returns zero to ten decimal places, but the same command
with $1$ replaced by $3$ returns 51.72756346.

\para{A rank six example.} More typically, {\em Magma} does not know $F_p(M,x)$ 
and $c_p$ for wild primes $p$ and one needs to input this information.
As an example, take 
 $M=H(\Phi_3^4/\Phi_1^8,1)$ with Hodge
 vector $(1,1,1,0,0,1,1,1)$.
 The only prime bad for the data is $p=3$.  
 A good first guess is that $F_3(M,x)$ is just the 
 constant $1$.      A short search over 
 some possible $c_3$ is implemented after redefining \verb@Q@ by
 \smallskip
 
 \verb@[CFENew(LSeries(Q,1:Precision:=10),@
 
 \verb@ BadPrimes:=[<3,c,1>]): c in [6..10]];@
 \smallskip
 
 \noindent The returned number for $c=9$ is $0.0000000000$, 
 while the numbers for the other $c$ are all at least $0.1$.  
 This information strongly suggests that indeed $F_3(M,x)=1$
 and $c_3=9$.   After setting up \verb@L@ with
 \verb@[<3,9,1>]@, analytic calculations can 
 be done as before.   For example, here 
 the order of central vanishing is apparently 2.  
 In the miraculous command \verb@CFENew@,
 \verb@CFE@ stands for the {\em Magma} command 
 \verb@CheckFunctionalEquation@, implemented
 by Tim Dokchitser using his \cite{Dok}; \verb@New@
 reflects subsequent improvements by Watkins.

\section{Bad primes}  
\label{bad} 
Fix a hypergeometric motive $M=H(\fa,t)$ and a prime $p$.  We now
sketch how {\em Magma} computes the local invariants 
when $p$ is tame for $(\fa,t)$, and describe some 
conjectural basic features for the case when $p$ is wild 
for $(\fa,t)$.  

\para{Tame primes.}  
When $p$ is tame for $(\fa,t)$, the conductor exponent $c_p$ is 
the codimension of the invariants of a power of a Levelt matrix $h_\tau$
from Section~\ref{monodromy}.  When $\ord_p(t-1) \geq 1$,
the simple shape 
of $h_1$ gives a completely explicit formula: $c_p=1$
except in the orthogonal case with $\ord_p(t-1)$ even, 
where $c_p=0$.  When $\ord_p(t-1) \leq 0$, 
\begin{equation}
c_p =  \mbox{rank}(h_\tau^{|k|} -I).  \label{cptame}
\end{equation}
Here $k=\ord_p(t)$, $\tau = \infty$ if $k$ is negative, and $\tau=0$ if 
$k$ is positive.   So there is separate periodic behavior for
$k<0$ and $k>0$, as illustrated by the top part of Figure~\ref{pict53both}.
The example of this table comes from the case $(a,b)=(3,5)$ of \eqref{123},
so the conductor there is very simply computed as the discriminant
of the octic algebra $\Q[x]/(5 x^8 + 8 t x^5 + 3 t^3)$.  

Because ramification is at worst tame, the degree of $F_p(M,x)$ 
is $n-c_p$.   When $\ord_p(t-1)$ is
positive, $F_p(M,x)$ is computed by slightly modifying
the formulas for point counts sketched in Section~\ref{sect:arithmetic}.
In the other cases, $F_p(M,x)$ comes
from Jacobi motives as mentioned around \eqref{fermat}, 
extracted from how 
the family $X_{\fa,t}$ degenerates at the relevant
cusp $\tau \in \{0,\infty\}$.

\para{Wild primes.} To simplify the overview, we just exclude the case 
where $\ord_p(t-1) \geq 1$.  Write specialization points as 
$t = vp^k$ with $k = \ord_p(t)$.  The bottom part of Figure~\ref{pict53both} shows
right away that the situation is complicated.   
\begin{figure}[htb]
\begin{center}
\includegraphics[width=2.8in]{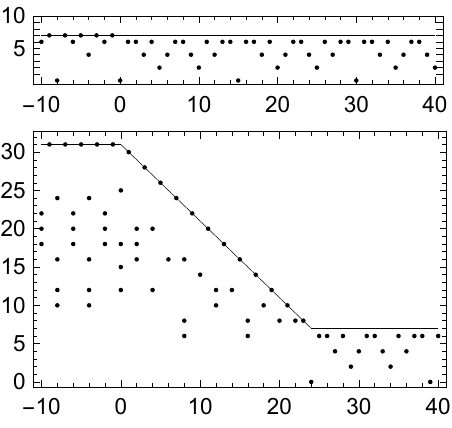}
\end{center}
\caption{\label{pict53both} Pairs $(k,c_p)$ where $k = \ord_p(t)$
and $c_p = \ord_p(\mbox{Conductor}(H([-8,3,5],t)))$,  compared with the graph
of the corresponding $\sigma$.    
 Top: The tame cases $p>5$.
Bottom: The wild case $p=2$.} 
\end{figure}

A function $\sigma$ is graphed in both parts of Figure~\ref{pict53both} 
and its general definition goes as follows.  
For $d$ a positive integer, write
\[
s(d) =
\left\{
\begin{array}{ll}
1 & \mbox{if $\gcd(d,p)=1$,} \\
1+\ord_p(d) + \frac{1}{p-1}  & \mbox{else.}
\end{array}
\right.
\]
Let
\begin{align*}
\sigma_\infty & = \sum_{i=1}^n s(\mbox{denom}(\alpha_i)), &
\sigma_0 & = \sum_{i=1}^n s(\mbox{denom}(\beta_i)). 
\end{align*}
Define $k_{\rm crit} = \sigma_\infty-\sigma_0 = -\sum_j \gamma_j \ord_p(\gamma_j)$ 
and transition points $k_\infty = \min(k_{\rm crit},0)$ 
and $k_0 = \max(k_{\rm crit},0)$.  Then 
\[
\sigma(k) = 
\left\{
\begin{array}{ll}
\sigma_\infty & \mbox{if $k \leq k_\infty$},\\
\max(\sigma_\infty,\sigma_0)-|k|& \mbox{if $k_\infty \leq k \leq k_0$}, \\
\sigma_0 & \mbox{if $k \geq k_0$}. \\
\end{array}
\right.
\]
In the tame case, 
$\sigma$ is just the constant function $n$. 
In general, there are plateaus corresponding
to the cusps $\infty$ and $0$, and then
a ramp of length $|k_{\rm crit}|$ between them.  

We conjecture that 
\begin{equation}
\label{rampbound}
c_p \leq \sigma(k)-\mbox{degree}(F_p(M,x)),
\end{equation}
with equality if
$k$ and $p$ are relatively prime.   
The second statement is proved in \cite{LNV}
in the general trinomial setting of 
\eqref{123}.  All of \eqref{rampbound}
has been computationally verified 
in many instances.  As one passes from 
one family to another via mod $\ell$
congruences as in Section~\ref{frobenius}, wild ramification
at $p$ does not change.  This fact and other
theoretical stabilities give us confidence in 
\eqref{rampbound}.   To make
{\em Magma} more fully automatic, a key step
would be to define a more 
 complicated
function $\sigma(k,v)$ with
$c_p \leq \sigma(k,v)-\mbox{degree}(F_p(M,x))$, and equality
under broad circumstances. 

 At present, we understand a factor $f_p(M,x)$ 
of the Frobenius polynomial $F_p(M,x)$ as
follows.  For $k \neq k_{\rm crit}$, $f_p(M,x)$ 
comes from modifying the tame formulas; 
in particular its degree is given by replacing
 $k$ by $k-k_{\rm crit}$ in \eqref{cptame}. 
 If $k=k_{\rm crit}$, corresponding to being
 at the bottom of the ramp, we use an {\em erasing} principle
 explained to us by Katz.  
Here one simply ignores all $\alpha_j$ and
$\beta_j$ that have denominator divisible by $p$.   Let
$n_\infty$ and $n_0$ be respectively the number of 
$\alpha_j$'s and $\beta_j$'s remaining.  Then
$n_\infty-n_0$ is a multiple of $p-1$, so that the formulas 
described in Section~\ref{sect:arithmetic} still make sense, as 
the choice of an auxiliary additive character on $\F_p$
again does not matter.  The resulting $f_p(M,x)$ has
degree $\max(n_\infty,n_0)$.  
We conjecture that the complementary factor 
$F_p(M,x)/f_p(M,x)$ is $1$ whenever $p$ and $k$ are relatively 
prime.  In practice, when $p \mid k$ it is usually $1$ also, but 
not always.  

\section{Automorphy} 
\label{automorphy}  
One of the most exciting aspects of the theory of motives is its conjectured
extremely tight connection to automorphic representations of 
adelic groups through the Langlands program.   

\para{Background.}  Let $\A$ be the adele ring of $\Q$; it is a restricted product of all
the completions $\Q_p$, including $\Q_\infty=\R$.  A cuspidal automorphic representation of 
$GL_n(\A)$ has an $L$-function known to have an analytic continuation and functional equation. 
The main conjecture is that,
after incorporating Tate twists to make normalizations match, the set 
of $L$-functions coming from irreducible rank $n$ motives in 
$\cM(\Q,\C)$ is exactly the subset of automorphic $L$-functions for
which the infinity factor has the form \eqref{inffact}.   

\para{The case $n=2$.}  For a motive $M \in \cM(\Q,\C)$ with 
nonvanishing Hodge numbers $h^{w,0}=h^{0,w}=1$ and conductor
$N$, one can switch to classical language.  The desired automorphic 
representation is entirely given by a power series in $q=e^{2 \pi i z}$ as in \eqref{legendre-mod},
but now this newform on $\Gamma_0(N)$ has weight $w+1$.

To exhibit some matches between motivic and automorphic $L$-functions,
consider the four reflexive parameters $\fa$ yielding motives $H(\fa,1)$
with Hodge vector $(1,1,0,0,1,1)$:
 \[
{\renewcommand{\arraycolsep}{3pt}
\begin{array}{c|rrrr|rrrr}
\fa & a_5 & a_7 & N' && b_5 & b_7 & N'' &\\
\hline
\Phi_2^6/\Phi_1^3 & -2 & 24 & 8 &\eta_2^4 \eta_4^4 & 74 & 24 & 8 &\\
\Phi_2^4\Phi_3 / \Phi_1^4\Phi_6 & -18 & 8 & 12 && 54 & -88 & 4 & \eta_2^{12} \\
\Phi_2^2\Phi_3^2/\Phi_1^2\Phi_6^2 & -6 & -16 & 18 && -66 & 176 & 6 &  \\
\Phi_3^3/\Phi_6^3 & -16 & -12 & 72 && -16 & 12 & 72 &
\end{array}
}
\]
{\em Magma} computes automatically with these reducible motives, reporting their 
conductors to be $N=64$, $48$, $108$, and $5184$.  
However these computations do not see the 
decompositions $H(\fa,1)=M'(-1) \oplus M''$ analogous to \eqref{dyadic16},
where now 
$M'$ and $M''$ respectively have Hodge vectors $(1,0,0,1)$ and
$(1,0,0,0,0,1)$.
In the Frobenius polynomial 
\[
F_p(M(\fa,1),x) = (1 - p a_p x + p^5 x^2) (1 - b_p x + p^5 x^2),
\] 
the $pa_p$ belonging to $M'(-1)$ can be distinguished
from the $b_p$ belonging to $M''$ 
whenever the latter is not a multiple of $p$.   The reader
might enjoy searching in the LMFDB's complete lists \cite{lmfdb}
of modular forms to see that the $a_p$ and $b_p$ for $p=5$ and $p=7$
let one identify the relevant forms and 
in particular determine
the above-displayed factorizations $N=N' N''$.  
Part of the further information given by the LMFDB is that
 two of the forms are expressible using the Dedekind eta function $\eta_1$,
via $\eta_d = q^{d/24} \prod_{j=1}^{\infty} (1-q^{dj})$.  

\para{Higher rank.} For a given motive $M \in \cM(\Q,\C)$ with
larger rank $n$, one can usually replace $GL_n(\A)$ by the adelic points
of a smaller group determined by the motivic Galois group
$G$ of $M$.  In favorable cases, the representation sought
again corresponds to a holomorphic form.  For rank three orthogonal motives, classical modular
forms are again relevant, but a symmetric square is now involved.  In rank four, 
Hilbert modular forms are needed for orthogonal motives and Siegel
modular forms are needed for symplectic motives.  
Numerical and sometimes proved matches have been found in these three settings.
For example, \cite{DPVZ} treats some interesting rank four orthogonal cases.

Generally speaking, the Hodge numbers 
of central concern earlier in this survey continue to play a large role.
In particular, motives for which all $h^{p,q}$ are $0$ or $1$ have theoretical
advantages, and their
motivic $L$-functions at least have a meromorphic
continuation with the right functional equation \cite{PT}.

\section{Numerical computations}  
\label{numerical}
We promised in the introduction that we would equip  
the reader to numerically explore a large collection of motivic $L$-functions.
We conclude this survey by giving sample computations in the context of 
two important topics, always assuming that the expected analytic continuation 
and functional equation indeed hold.    In both topics, we let $c = \frac{1}{2}+\frac{w}{2}$ be the center of the functional equation.
The conductors $N$ in our examples are 
small for their Hodge vectors $h$, allowing us to keep runtimes short and/or
work to high precision.

\para{Special values.} 
If $M$ is a motive in $\cM(\Q,\Q)$ then the numbers $L(M,k)$ for integers $k \leq c$ are 
mostly forced to be $0$, because of poles in the infinity factor \eqref{inffact}
and the functional equation.  However, when $L(M,k)$ is nonzero it is 
expected to be arithmetically significant \cite{Del-PL}.  The arithmetic interpretation involves a 
determinant of periods like \eqref{euler-integral}.  To see the significance
without entering into periods, one can look at the ratio 
$r_d = L(M \otimes \chi_d,k)/L(M,k)$, for $d$ a positive quadratic
discriminant.  Then the periods cancel out so that
$r_d$ should be rational.

For a sample computation, take $M=H(\Phi_2^5/\Phi_1^5,2^{10})$ and use
\eqref{defineq} and  
\verb@L:=LSeries(Q,1024)@ to define its $L$-function, as usual.    
While $2$ is wild for the family, it is unramified in 
$M$ because because the exponent $10$ is at the bottom 
of the ramp of Section~\ref{bad}.  The erasing procedure
from the end of that section applies, yielding  
\[
F_2(M,x) = (1-4x)(1+5x+10x^2+80 x^3 + 256 x^4).
\] 
Since 
$t - 1 =  1023 =  3 \cdot 11 \cdot 31$ is squarefree, it is the conductor,  by the recipe
before \eqref{cptame}.  {\em Magma} gets all the bad factors right automatically.  As a confirmation,
 \verb@CFENew(L)@ quickly returns $0$ to the default
$30$ digits.

\verb@Evaluate(L,2)@ gives
$0.4278180899 \cdots$.
Twisting by a $d$ with $\gcd(d,1023)=1$ makes the conductor go up by a factor 
of $d^5$ and precision needs to be reduced.  
\smallskip

\verb@Evaluate(LSeries(Q,1024:@

\verb@QuadraticTwist:=5,Precision:=10),2);@
\smallskip

\noindent takes six minutes to give its answer of $35.04685793$.  This ratio
and then two others are apparently
\begin{align*}
r_5 & = \frac{2^{11}}{5^2}, & r_8 & = {2^6 \, 5}, & 
r_{13} &
=\frac{2^{10} \,251}{13^2} .
\end{align*}
The two $L$-functions appearing in $r_d$ are completely different
analytically, and so the apparent fact that quotients are rational
is very remarkable.  

Readers wanting to work out their own examples might want
to begin with $M$ having odd weight.  Then if $L(M,c) \neq 0$,
one has conjecturally rational quotients $r_d$ for $k=c$.
The lateral argument $k=c-j$ fits into the theory
only in the rare case that the $2j$ most central entries of the Hodge vector are $0$.    
In the even weight case, one needs to have $h_+=0$ to 
make $k=c-\frac{1}{2}$ fit into the theory,  as in our example.

\para{Critical zeros.}  For a 
weight $w$ 
motive $M$, all the zeros of the 
completed $L$-function $\Lambda(M,s)$ lie in the critical strip
$c- \frac{1}{2} \leq \Real(s) \leq c+\frac{1}{2}$. 
 The Riemann hypothesis
for $M$ then predicts that all the zeros lie on the critical line
$\Real(s) = c$.  We now show by examples that numerical identification of low-lying zeros is 
possible in modestly high rank.  

For the examples, take $M_{10,w} = H(\fa_w,1)$ with 
$\fa_{10} = {\Phi_4 \Phi_2^9}/{\Phi^{11}_1}$ and 
$\fa_7  = {\Phi_4^4 \Phi_2^4}/{\Phi_8^2 \Phi_1^4}$.
So $M_{10,10}$ is orthogonal with Hodge vector $(1,1,1,1,1,0,1,1,1,1,1)$
while $M_{10,7}$ is symplectic with Hodge vector $(1,1,2,1,1,2,1,1)$.

The only bad prime in each case is $2$.  A search says that $F_2(M_{10,10},x) = 1+32 x$ and
$c_2 = 11$. 
For $M_{10,7}$, $k=k_{\rm crit}=0$ so erasing
applies, yielding $1 + 4 x + 96 x^2 + 512 x^3 + 16384 x^4$
as a factor of $F_2(M_{10,7},x)$.  A short search says that this factor is all of 
$F_2(M_{10,7},x)$ and $c_2=18$.  

In general, the Hardy Z-function of a motive $M$ is
\[
Z(M,t) = \epsilon^{1/2} \frac{N^{s/2} \Gamma_{h,\sigma}(s)}{|N^{s/2} \Gamma_{h,\sigma}(s)|} L(M,s),
\]
with $s=c+it$.  It is a real-valued function of the real variable $t$, even or odd depending
on whether the sign $\epsilon$ is $1$ or $-1$.

 \begin{figure}[htb]
 \begin{center}
\includegraphics[width=3in]{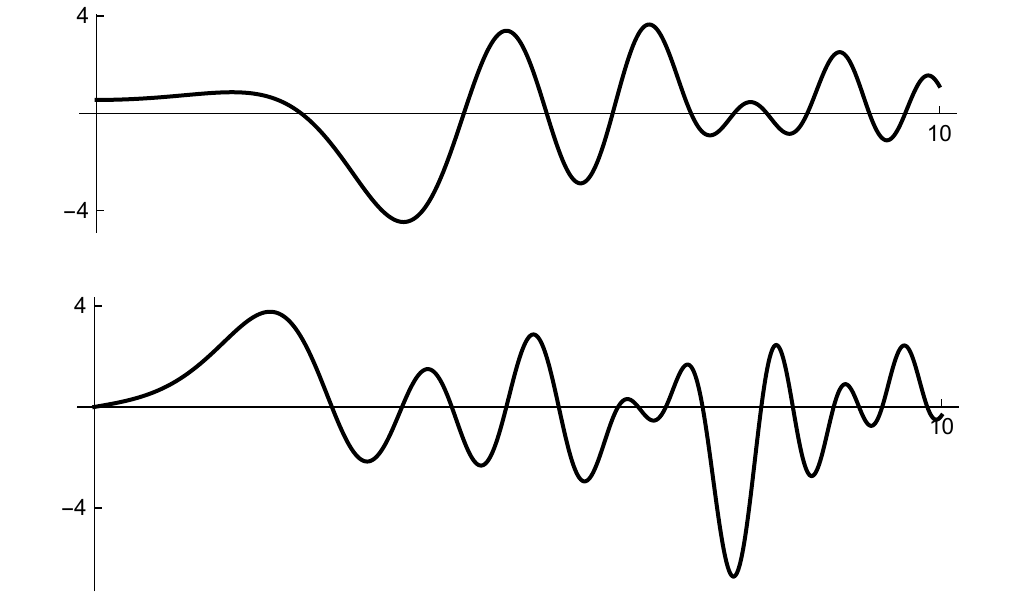}
\end{center}
\caption{\label{twoz} Graphs of $Z(M_{10,10},t)$ and $Z(M_{10,7},t)$}
\end{figure}

Figure~\ref{twoz} was computed via many calls to \verb@Evaluate@ at points of the 
form $c + i t$.   The signs in the two cases are $1$ and $-1$, and the orders of 
central vanishing are the minimum possible, $0$ and $1$.  
On both plots, all local maxima are above the axis and all local minima are
beneath the axis.  Zeros off the critical line would likely cause a disruption
of this pattern; thus the plots not only identify zeros on the critical line, but suggest
a lack of zeros off the critical line.
 \bibliography{HGMs}

\begin{bibdiv}
\begin{biblist}

\bib{And}{book}{
      author={Andr\'{e}, Yves},
       title={Une introduction aux motifs (motifs purs, motifs mixtes,
  p\'{e}riodes)},
      series={Panoramas et Synth\`eses},
   publisher={Soci\'{e}t\'{e} Math\'{e}matique de France, Paris},
        date={2004},
      volume={17},
        ISBN={2-85629-164-3},
      review={\MR{2115000}},
}

\bib{BCM}{article}{
      author={Beukers, Frits},
      author={Cohen, Henri},
      author={Mellit, Anton},
       title={Finite hypergeometric functions},
        date={2015},
        ISSN={1558-8599},
     journal={Pure Appl. Math. Q.},
      volume={11},
      number={4},
       pages={559\ndash 589},
  url={https://doi-org.ezproxy.morris.umn.edu/10.4310/PAMQ.2015.v11.n4.a2},
      review={\MR{3613122}},
}

\bib{BH}{article}{
      author={Beukers, F.},
      author={Heckman, G.},
       title={Monodromy for the hypergeometric function {$_nF_{n-1}$}},
        date={1989},
        ISSN={0020-9910},
     journal={Invent. Math.},
      volume={95},
      number={2},
       pages={325\ndash 354},
         url={https://doi-org.ezproxy.morris.umn.edu/10.1007/BF01393900},
      review={\MR{974906}},
}

\bib{CG}{article}{
      author={Corti, Alessio},
      author={Golyshev, Vasily},
       title={Hypergeometric equations and weighted projective spaces},
        date={2011},
        ISSN={1674-7283},
     journal={Sci. China Math.},
      volume={54},
      number={8},
       pages={1577\ndash 1590},
  url={https://doi-org.ezproxy.morris.umn.edu/10.1007/s11425-011-4218-5},
      review={\MR{2824960}},
}

\bib{Del-PL}{inproceedings}{
      author={Deligne, P.},
       title={Valeurs de fonctions {$L$} et p\'{e}riodes d'int\'{e}grales},
        date={1979},
   booktitle={Automorphic forms, representations and {$L$}-functions ({P}roc.
  {S}ympos. {P}ure {M}ath., {O}regon {S}tate {U}niv., {C}orvallis, {O}re.,
  1977), {P}art 2},
      series={Proc. Sympos. Pure Math., XXXIII},
   publisher={Amer. Math. Soc., Providence, R.I.},
       pages={313\ndash 346},
      review={\MR{546622}},
}

\bib{Dok}{article}{
      author={Dokchitser, Tim},
       title={Computing special values of motivic {$L$}-functions},
        date={2004},
        ISSN={1058-6458},
     journal={Experiment. Math.},
      volume={13},
      number={2},
       pages={137\ndash 149},
  url={http://projecteuclid.org.ezproxy.morris.umn.edu/euclid.em/1090350929},
      review={\MR{2068888}},
}

\bib{DPVZ}{article}{
      author={Demb\'el\'e, Lassina},
      author={Panchishkin, Alexei},
      author={Voight, John},
      author={Zudilin, Wadim},
       title={Special hypergeometric motives and their {L}-functions: Asai
  recognition},
        date={2020},
     journal={Experimental Mathematics},
      volume={0},
      number={0},
       pages={1\ndash 23},
      eprint={https://doi.org/10.1080/10586458.2020.1737990},
         url={https://doi.org/10.1080/10586458.2020.1737990},
}

\bib{Fed}{article}{
      author={Fedorov, Roman},
       title={Variations of {H}odge structures for hypergeometric differential
  operators and parabolic {H}iggs bundles},
        date={2018},
        ISSN={1073-7928},
     journal={Int. Math. Res. Not. IMRN},
      number={18},
       pages={5583\ndash 5608},
         url={https://doi-org.ezproxy.morris.umn.edu/10.1093/imrn/rnx044},
      review={\MR{3862114}},
}

\bib{GKZ}{book}{
      author={Gelfand, I.~M.},
      author={Kapranov, M.~M.},
      author={Zelevinsky, A.~V.},
       title={Discriminants, resultants, and multidimensional determinants},
   publisher={Birkh\"{a}user Boston, Inc., Boston, MA},
        date={1994},
        ISBN={0-8176-3660-9},
  url={https://doi-org.ezproxy.morris.umn.edu/10.1007/978-0-8176-4771-1},
      review={\MR{1264417}},
}

\bib{RV-bez}{article}{
      author={Jouve, F.},
      author={Rodriguez~Villegas, F.},
       title={On the bilinear structure associated to {B}ezoutians},
        date={2014},
        ISSN={0021-8693},
     journal={J. Algebra},
      volume={400},
       pages={161\ndash 184},
  url={https://doi-org.ezproxy.morris.umn.edu/10.1016/j.jalgebra.2013.10.030},
      review={\MR{3147370}},
}

\bib{Katz-ESDE}{book}{
      author={Katz, Nicholas~M.},
       title={Exponential sums and differential equations},
      series={Annals of Mathematics Studies},
   publisher={Princeton University Press, Princeton, NJ},
        date={1990},
      volume={124},
        ISBN={0-691-08598-6; 0-691-08599-4},
         url={https://doi-org.ezproxy.morris.umn.edu/10.1515/9781400882434},
      review={\MR{1081536}},
}

\bib{Katz-RLS}{book}{
      author={Katz, Nicholas~M.},
       title={Rigid local systems},
      series={Annals of Mathematics Studies},
   publisher={Princeton University Press, Princeton, NJ},
        date={1996},
      volume={139},
        ISBN={0-691-01118-4},
         url={https://doi-org.ezproxy.morris.umn.edu/10.1515/9781400882595},
      review={\MR{1366651}},
}

\bib{lmfdb}{misc}{
      author={{LMFDB Collaboration}, The},
       title={The {L}-functions and modular forms database},
         how={\url{http://www.lmfdb.org}},
        date={2021},
        note={\url{http://www.lmfdb.org}},
}

\bib{LNV}{article}{
      author={Llorente, P.},
      author={Nart, E.},
      author={Vila, N.},
       title={Discriminants of number fields defined by trinomials},
        date={1984},
        ISSN={0065-1036},
     journal={Acta Arith.},
      volume={43},
      number={4},
       pages={367\ndash 373},
         url={https://doi-org.ezproxy.morris.umn.edu/10.4064/aa-43-4-367-373},
      review={\MR{756288}},
}

\bib{PT}{article}{
      author={Patrikis, Stefan},
      author={Taylor, Richard},
       title={Automorphy and irreducibility of some {$l$}-adic
  representations},
        date={2015},
        ISSN={0010-437X},
     journal={Compos. Math.},
      volume={151},
      number={2},
       pages={207\ndash 229},
  url={https://doi-org.ezproxy.morris.umn.edu/10.1112/S0010437X14007519},
      review={\MR{3314824}},
}

\bib{Rob-polys}{misc}{
      author={Roberts, David~P.},
       title={Shioda polynomials for $\!\! {W}({E}_n) \!\!$ {B}eukers-{H}eckman
  covers},
        date={2018},
         url={https://www.davidproberts.net/presentations},
        note={\url{http://www.davidproberts.net}},
}

\bib{Rob}{incollection}{
      author={Roberts, David~P.},
       title={The explicit formula and a motivic splitting},
        date={2019},
   booktitle={Matrix {A}nnals},
      series={MATRIX Book Series},
      volume={2},
   publisher={Springer, Cham},
       pages={429\ndash 433},
         url={xxx},
}

\bib{RV}{incollection}{
      author={Rodriguez~Villegas, Fernando},
       title={Hypergeometric families of {C}alabi-{Y}au manifolds},
        date={2003},
   booktitle={Calabi-{Y}au varieties and mirror symmetry ({T}oronto, {ON},
  2001)},
      series={Fields Inst. Commun.},
      volume={38},
   publisher={Amer. Math. Soc., Providence, RI},
       pages={223\ndash 231},
      review={\MR{2019156}},
}

\bib{RV-Mixed}{article}{
      author={Rodriguez~Villegas, Fernando},
       title={Mixed {H}odge numbers and factorial ratios},
        date={2019},
     journal={arxiv 1907.02722},
}

\bib{Wat}{misc}{
      author={Watkins, Mark},
       title={Hypergeometric motives package},
        date={2015},
        note={Available in Magma, including through the online calculator at
  \url{http://magma.maths.usyd.edu.au/calc/}},
}

\end{biblist}
\end{bibdiv}
\end{document}